\documentclass[12pt,twoside]{article}

\usepackage{amssymb,amsmath,amsfonts,graphicx,makeidx}

\allowdisplaybreaks

\newenvironment{Proof}{\removelastskip\par\medskip
\noindent{\em Proof.} \rm}{\penalty-20\null\hfill$\square$\par\medbreak}

\newcommand{\ee}{\mathbb{E}}
\newcommand{\Dom}{{\rm Dom \ \!}} 
\newtheorem{assumption}{Assumption}[section]{\bfseries}{\rmfamily}
\newcommand{\ignore}[1]{}

\def\real{{\mathord{\mathbb R}}}
\def\inte{{\mathord{\mathbb N}}}

\def\Ent{{\mathrm{{\rm Ent \ \! }}}}
\def\PP{{\mathord{{\rm I\kern-2.8pt P}}}}

\numberwithin{equation}{section}

\newtheorem{prop}[equation]{Proposition}
\newtheorem{lemma}[equation]{Lemma}
\newtheorem{definition}[equation]{Definition}
\newtheorem{corollary}[equation]{Corollary}
\newtheorem{theorem}[equation]{Theorem}
\newtheorem{remark}[equation]{Remark}

\def\Dom{{\mathrm{{\rm Dom  \ \!}}}}

\def\Ent{{\mathrm{{\rm Ent \ \!}}}}
\def\var{{\mathrm{{\rm var \ \!}}}}

\def\Cov{{\mathrm{{\rm Cov \ \!}}}}
\def\var{{\mathrm{{\rm var \ \!}}}}
\def\Var{{\mathrm{{\rm Var \ \!}}}}
\def\Vect{{\mathrm{{\rm Vect \ \!}}}}

\def\d{{\mathord{{\rm d}}}}

\def\P{\mathbb{P}}

\def\E{\mathop{\hbox{\rm I\kern-0.20em E}}\nolimits}

\def\Var{\mathop{\hbox{\rm Var}}\nolimits}
\def\Cov{\mathop{\hbox{\rm Cov}}\nolimits}

\title{
  \huge
  Stochastic analysis of Bernoulli processes 
}

\author{
\large 
Nicolas Privault 
\\ 
\normalsize 
Division of Mathematical Sciences 
\\ 
\normalsize 
School of Physical and Mathematical Sciences 
\\ 
\normalsize 
Nanyang Technological University 
\\ 
\normalsize 
21 Nanyang Link 
\\ 
\normalsize 
Singapore 637371
\\ 
\normalsize 
nprivault@ntu.edu.sg
}

\begin{document}

\maketitle

\begin{abstract} 
 These notes survey some aspects of discrete-time chaotic calculus 
 and its applications, based on the chaos representation property 
 for i.i.d. sequences of random variables. 
 The topics covered include the Clark formula and predictable 
 representation, anticipating calculus, covariance identities and 
 functional inequalities (such as deviation and logarithmic Sobolev 
 inequalities), and an application to option hedging in discrete time. 
\end{abstract} 
\\ 

\noindent 
 Classification: 60G42, 60G50, 60G51, 60H30, 60H07. 

\baselineskip0.7cm

\section{Introduction} 

 Stochastic analysis can be viewed as an infinite-dimensional 
 version of classical analysis, developed in relation to stochastic 
 processes. 

 In this survey we present a construction of the basic operators 
 of stochastic analysis (gradient and divergence) 
 in discrete time for Bernoulli processes. 
 Our presentation is based on the chaos representation property 
 and discrete multiple stochastic integrals with respect to 
 i.i.d. sequences of random variables. 
 The main applications presented are to functional inequalities 
 (deviation inequalities, logarithmic Sobolev inequalities) 
 in discrete settings, 
 cf. \cite{gaopri}, \cite{hp}, \cite{prisch2}, and to option 
 pricing and hedging in discrete time mathematical finance. 

 Other approaches to discrete-time stochastic analysis 
can be found in Holden et al. \cite{holden1} (1992), 
\cite{holden2} (1993), Leitz-Martini \cite{leitz} (2000), 
and also in Attal \cite{attal1} (2003) in the framework of 
quantum stochastic calculus, see also the recent paper \cite{gzyl} 
 by H. Gzyl (2005). 
 
 This survey can be roughly divided into a first part (Sections~\ref{s2} to 
 \ref{s8}) in which we present the main basic results 
 and analytic tools, and a second part (Sections~\ref{devsec} to \ref{hdg}) 
 which is devoted to applications.   

 We proceed as follows. 
 In Section~\ref{s2} we consider a family of discrete-time 
 normal martingales. 
 The next section is devoted to the construction of 
 the stochastic integral of predictable square-integrable 
 processes with respect to such martingales. 
 In Section~\ref{s3} we construct the associated multiple 
 stochastic integrals of symmetric functions on $\inte^n$, 
 $n\geq 1$. 
 Starting with Section~\ref{s3.1} we focus on a particular 
 class of normal martingales satisfying a structure equation. 
 The chaos representation property is studied in Section~\ref{s3.2} 
 in the case of discrete time random walks with independent 
 increments. 
 A gradient operator $D$ acting by finite differences is 
 introduced in Section~\ref{s4} in connection with multiple 
 stochastic integrals, and used in Section~\ref{s5} 
 to state a Clark predictable representation formula. 
 The divergence operator $\delta$, adjoint of $D$, is 
 presented in Section~\ref{s6} as an extension of the 
 discrete-time stochastic integral. 
 It is also used in Section~\ref{s7} to express 
 the generator of the Ornstein-Uhlenbeck process. 
 Covariance identities are stated in Section~\ref{s8}, 
 both from the Clark representation formula and by 
 use of the Ornstein-Uhlenbeck semigroup. 

 Functional inequalities on Bernoulli space are presented 
 as an application in Sections~\ref{devsec} and \ref{lsid}. 
 On the one hand, in Section~\ref{devsec} 
 we prove several deviation inequalities 
 for functionals of an infinite number of i.i.d. Bernoulli 
 random variables. 
 Then in Section~\ref{lsid} we state different versions of 
 the logarithmic Sobolev inequality in discrete settings 
 (modified, $L^1$, sharp) which allow one to control the 
 entropy of random variables. 
 In particular we recover and extend some results 
 of \cite{bobkov}, using the method of \cite{gaopri}. 
 Our approach is based on the intrinsic tools 
 (gradient, divergence, Laplacian) of 
 infinite-dimensional stochastic analysis. 
 We refer to \cite{bht2}, \cite{bobkovsyk}, \cite{ht2},
 \cite{ledouxesaim}, for other versions of logarithmic 
 Sobolev inequalities in discrete settings, 
 and to \cite{daipra}, \cite{wuls2} for the Poisson case. 

 Section~\ref{s11} contains a change of variable formula 
 in discrete time, which is applied 
 with the Clark formula in Section~\ref{hdg} 
 to a derivation of the Black-Scholes formula 
 in discrete time, i.e. in the Cox-Ross-Rubinstein model, see e.g. 
 \cite{lamberton}, $\S$15-1 of \cite{williams}, or \cite{ruiz}, 
 for other approaches. 

\section{Discrete-Time Normal Martingales}\index{normal martingale!discrete time} 
\label{s2} 
 Consider a sequence $(Y_k)_{k\in \inte}$ 
 of (not necessarily independent) random 
 variables on a probability space $(\Omega , {\cal F} , \P)$. 
 Let $({\cal F}_n)_{n\geq -1}$ denote the filtration generated by 
 $(Y_n)_{n\in \inte}$, i.e. 
$$ 
 {\cal F}_{-1}=\{\emptyset , \Omega \} 
, 
$$ 
 and 
$${\cal F}_n = \sigma (Y_0,\ldots ,Y_n) 
, \qquad n\geq 0
. 
$$ 
 Recall that a random variable $F$ is said to be 
 ${\cal F}_n$-measurable if it can be written as a function 
$$F=f_n (Y_0,\ldots,Y_n)$$ 
 of $Y_0,\ldots , Y_n$, where 
 $f_n : \real^{n+1} \to \real$. 
\bigskip
 
\begin{assumption} 
 We make the following assumptions on the sequence 
 $(Y_n)_{n\in\inte}$: 
\begin{description} 
\item{a)} 
 it is conditionally centered: 
\begin{equation} 
\label{dse0} 
\ee [Y_n \mid {\cal F}_{n-1}]=0, \qquad 
n\geq 0, 
\end{equation} 
\item{b)} its conditional quadratic variation satisfies: 
$$\ee[Y_n^2 \mid {\cal F}_{n-1} ]=1, \qquad 
 n\geq 0. 
$$ 
\end{description} 
\end{assumption} 
 Condition \eqref{dse0} implies that the process 
 $(Y_0+\cdots +Y_n)_{n\geq 0}$ is an ${\cal F}_n$-martingale. 
 More precisely, the sequence $(Y_n)_{n\in\inte}$ and the process 
 $(Y_0+\cdots +Y_n)_{n\geq 0}$ can be viewed 
 respectively as a (correlated) noise and as a normal martingale 
 in discrete time. 
\section{Discrete Stochastic Integrals}\index{stochastic integral!discrete time}
 In this section we construct the discrete stochastic integral 
 of predictable square-summable processes with respect to a 
 discrete-time normal martingale. 
\begin{definition} 
 Let $(u_k)_{k\in \inte}$ be a uniformly bounded sequence of 
 random variables with finite support in $\inte$, i.e. there 
 exists $N\geq 0$ such that $u_k=0$ for all $k\geq N$. 
 The stochastic integral $J(u)$ of 
 $(u_n)_{n\in\inte}$ is defined as
$$J( u) = \sum_{k=0}^\infty u_kY_k 
. 
$$ 
\end{definition} 
 The next proposition states a version of 
 the It\^o isometry in discrete time. 
 A sequence $( u_n)_{n\in\inte}$ of random variables
 is said to be ${\cal F}_n$-predictable\index{predictable process} 
 if $u_n$ is ${\cal F}_{n-1}$-measurable
 for all $n\in \inte$, in particular $u_0$ is constant 
 in this case. 
\begin{prop} 
\label{is} 
 The stochastic integral operator $J(u)$ 
 extends to square-integrable predictable processes 
 $(u_n)_{n\in\inte} \in L^2 (\Omega \times \inte )$ 
 via the (conditional) isometry formula\index{It\^o!isometry}
\begin{equation} 
\label{iiiso} 
 \ee [ | J( {\bf 1}_{[n,\infty )} u) |^2 | \mid {\cal F}_{n-1} ] 
 = \ee[\Vert {\bf 1}_{[n,\infty )} u \Vert_{\ell^2(\inte)}^2 \mid {\cal F}_{n-1} ], 
 \qquad n\in\inte.
\end{equation}
\end{prop}
\begin{Proof} 
 Let $(u_n)_{n\in\inte}$ and $(v_n)_{n\in\inte}$ be bounded 
 predictable processes with finite support in $\inte$. 
 The product $u_kY_kv_l$, $0\leq k < l$, is ${\cal F}_{l-1}$-measurable, 
 and $u_kY_lv_l$ is ${\cal F}_{k-1}$-measurable, $0\leq l < k$. 
 Hence 
\begin{eqnarray*}
\lefteqn{ 
 \ee\left[ \sum_{k=n}^\infty u_kY_k  \sum_{l=0}^\infty v_lY_l
 \Big| {\cal F}_{n-1} \right] 
 = 
\ee\left[ \sum_{k,l=n}^\infty u_kY_k v_l Y_l 
\Big| {\cal F}_{n-1} \right]
} 
\\
& = &
\ee\left[ 
 \sum_{k=n}^\infty u_kv_k Y_k^2
 + \sum_{n\leq k<l } u_kY_k v_lY_l
 + \sum_{n\leq l < k} u_kY_k  v_lY_l
\Big| {\cal F}_{n-1} \right]
\\ 
& = &
 \sum_{k=n}^\infty \ee[ 
 \ee[u_kv_k Y_k^2 \mid {\cal F}_{k-1} ] 
 \mid {\cal F}_{n-1} ]
+ \sum_{n\leq k<l } 
 \ee[ 
 \ee[ 
 u_kY_k v_lY_l \mid {\cal F}_{l-1} ] 
 \mid {\cal F}_{n-1} 
 ]
\\ 
 & & + 
 \sum_{n\leq l < k} 
 \ee[ 
 \ee[ 
 u_kY_k  v_lY_l \mid {\cal F}_{k-1} ] 
 \mid {\cal F}_{n-1} 
 ]
\\ 
& = &
 \sum_{k=0}^\infty \ee[ 
 u_kv_k 
 \ee[ Y_k^2 \mid {\cal F}_{k-1} ] 
 \mid {\cal F}_{n-1} 
 ]
 + 
 2 \sum_{n\leq k<l } 
 \ee[ 
 u_kY_k v_l 
 \ee[ 
 Y_l \mid {\cal F}_{l-1} ] 
 \mid {\cal F}_{n-1} 
 ] 
\\
 & = & 
 \sum_{k=n}^\infty \ee[ u_kv_k \mid {\cal F}_{n-1} ] 
\\ 
 & = & 
 \ee\left[\sum_{k=n}^\infty u_kv_k \Big| {\cal F}_{n-1} \right]. 
\end{eqnarray*} 
 This proves the isometry property \eqref{iiiso} for $J$. 
 The extension to $L^2 (\Omega \times \inte )$ follows then 
 from a Cauchy sequence argument. 
 Consider a sequence of bounded predictable 
 processes with finite support converging to $u$ in 
 $L^2(\Omega \times \inte )$, for example the sequence 
 $(u^n)_{n\in\inte}$ defined as 
$$ 
 u^n 
 = 
 ( u^n_k )_{k\in\inte} 
 = 
 ( u_k {\bf 1}_{\{ 0 \leq k \leq n\}} {\bf 1}_{\{ |u_k|\leq n \}} )_{k\in\inte}, 
 \qquad 
 n\in\inte 
. 
$$ 
 Then the sequence $(J(u^n))_{n\in\inte}$ is Cauchy and 
 converges in $L^2(\Omega )$, hence we may define 
$$ 
 J(u) : = \lim_{k \to \infty} J(u^k ). 
$$ 
 From the isometry property \eqref{iiiso} applied with $n=0$, 
 the limit is clearly independent of the choice of the approximating sequence 
 $(u^k)_{k\in\inte}$. 
\end{Proof} 
 Note that by bilinearity, \eqref{iiiso} can also be written as 
$$ 
 \ee [ J( {\bf 1}_{[n,\infty )} u) J( {\bf 1}_{[n,\infty )} v )  | 
 {\cal F}_{n-1} ] 
 = 
 \ee[ \langle 
 {\bf 1}_{[n,\infty )} u 
 , 
 {\bf 1}_{[n,\infty )} v 
 \rangle_{\ell^2(\inte)} 
 \mid {\cal F}_{n-1} ], 
 \qquad n\in\inte, 
$$
 and that for $n=0$ we get 
\begin{equation} 
\label{lkoi} 
 \ee [ J( u) J( v ) ] 
 = 
 \ee[ \langle 
 u 
 , 
 v 
 \rangle_{\ell^2(\inte)} 
], 
\end{equation} 
 for all square-integrable predictable processes 
 $u = (u_k)_{k\in \inte}$ and $v = (v_k)_{k\in \inte}$. 
\begin{prop} 
\label{isomp15} 
 Let $(u_k)_{k\in\inte} \in L^2 (\Omega \times \inte )$ 
 be a predictable square-integrable process. 
 We have 
$$ 
 \ee[J (u ) \mid {\cal F}_k] = J (u{\bf 1}_{[0,k]}), \quad 
 k\in \inte 
. 
$$ 
\end{prop} 
\begin{Proof} 
 It is sufficient to note that 
\begin{eqnarray*} 
 \ee[J (u ) \mid {\cal F}_k] & = & 
 \ee\left[ \sum_{i=0}^k  u_iY_i \Big| {\cal F}_k \right] 
 + 
 \sum_{i=k+1}^\infty \ee\left[ u_iY_i \mid {\cal F}_k \right]  
\\ 
 & = & 
 \sum_{i=0}^k 
 u_iY_i 
 + 
 \sum_{i=k+1}^\infty 
 \ee\left[ \ee\left[ 
 u_iY_i \mid {\cal F}_{i-1} \right] \mid {\cal F}_k \right]  
\\ 
 & = & 
 \sum_{i=0}^k 
 u_iY_i 
 + 
 \sum_{i=k+1}^\infty 
 \ee\left[ 
 u_i 
 \ee\left[ 
 Y_i \mid {\cal F}_{i-1}  
 \right] \mid {\cal F}_k 
 \right]  
\\ 
 & = & 
 \sum_{i=0}^k 
 u_iY_i 
\\ 
 & = & 
 J (u{\bf 1}_{[0,k]}) 
. 
\end{eqnarray*} 
\end{Proof} 
\begin{corollary} 
 The indefinite stochastic integral 
 $(J (u{\bf 1}_{[0,k]}))_{k\in\inte}$ is a discrete time martingale with 
 respect to $({\cal F}_n)_{n\geq -1}$. 
\end{corollary} 
\begin{Proof} 
 We have 
\begin{eqnarray*} 
 \ee[J (u {\bf 1}_{[0,k+1]} ) \mid {\cal F}_k] 
 & = & 
 \ee[ \ee[ 
 J (u{\bf 1}_{[0, k + 1 ]}) \mid {\cal F}_{k+1} 
 \mid {\cal F}_k] 
\\ 
 & = & 
 \ee[ \ee[ 
 J (u ) \mid {\cal F}_{k+1} 
 \mid {\cal F}_k] 
\\ 
 & = & 
 \ee[ J (u ) \mid {\cal F}_k] 
\\ 
 & = & 
 J (u {\bf 1}_{[0,k]} )
. 
\end{eqnarray*} 
\end{Proof} 
\section{Discrete Multiple Stochastic Integrals}\index{multiple stochastic integrals!discrete time}
\label{s3} 
\noindent 
 The role of multiple stochastic integrals 
 in the orthogonal expansions of random variables 
 is similar to that of polynomials 
 in the series expansions of functions of a real variable. 
 In some cases, multiple stochastic integrals can be expressed 
 using polynomials, for example Krawtchouk polynomials in 
 the symmetric discrete case with $p_n=q_n=1/2$, $n\in\inte$, 
 see Relation~\eqref{krw} below. 
\begin{definition} 
 Let $\ell^2(\inte)^{\circ n}$ denote the subspace
 of $\ell^2(\inte)^{\otimes n} = \ell^2 (\inte^n)$ made of 
 functions $f_n$ that are symmetric in $n$ variables, i.e. such that
 for every permutation $\sigma$ of $\{1,\ldots , n\}$,
$$f_n (k_{\sigma (1)},\ldots ,k_{\sigma (n)})
 = f_n (k_1,\ldots ,k_n), \quad k_1,\ldots ,k_n\in\inte.
$$
\end{definition} 
 Given $f_1\in l^2 (\inte)$ we let 
$$J_1(f_1 ) = J(f_1 ) = \sum_{k=0}^\infty f_1 (k) Y_k 
. 
$$ 
 As a convention we identify $\ell^2(\inte^0)$ to $\real$ and let $J_0(f_0)=f_0$, 
 $f_0\in \real$. 
 Let 
$$
 \Delta_{n}=\{(k_1,\ldots ,k_n)\in \inte^n \ :
 \ k_i\not= k_j, \ 1\leq i < j \leq n \}, 
 \qquad n\geq 1.$$
 The following proposition gives the definition of multiple stochastic 
 integrals by iterated
 stochastic integration of predictable processes in the sense of Proposition~\ref{is}. 
\begin{prop} 
\label{prop1.18} 
 The multiple stochastic integral $J_n (f_n)$ of $f_n\in \ell^2(\inte)^{\circ n}$, 
 $n\geq 1$, is defined as
$$ 
J_n (f_n)
 = \sum_{(i_1,\ldots ,i_n)\in \Delta_n}
 f_n(i_1,\ldots ,i_n) Y_{i_1}\cdots Y_{i_n} 
.  
$$
 It satisfies the recurrence relation 
\begin{equation}
\label{recrel2} 
 J_{n}(f_{n}) = 
 n\sum_{k=1}^{\infty } 
 Y_k J_{n-1}(f_{n}( * ,k) {\bf 1}_{[0,k-1]^{n-1}}( * )) 
\end{equation} 
 and the isometry formula 
\begin{equation}
\label{isois}
\ee[J_n(f_n)J_m(g_m)]
 = \left\{ \begin{array}{ll}
 n! \langle {\bf 1}_{ \Delta_n} f_n , g_m \rangle_{\ell^2(\inte)^{\otimes n}}
 & \mbox{if } n=m,
\\
 0 & \mbox{if } n\not= m.
\end{array}
\right.
\end{equation}
\end{prop} 
\begin{Proof}
 Note that we have 
\begin{eqnarray} 
\nonumber 
 J_n (f_n)
 & = & n! \sum_{0\leq i_1 < \cdots < i_n} 
 f_n(i_1,\ldots ,i_n) Y_{i_1}\cdots Y_{i_n} 
\\ 
\label{eg1} 
 & = & 
 n! \sum_{i_n = 0 ~}^\infty
 \sum_{0\leq i_{n-1} < i_n}
 \cdots
 \sum_{0\leq i_1 < i_2}
 f_n(i_1,\ldots ,i_n) Y_{i_1}\cdots Y_{i_n} 
. 
\end{eqnarray} 
 Note that since $0\leq i_1 < i_2 < \cdots < i_n$ and 
 $0 \leq j_1 < j_2 < \cdots < j_n$ we have 
$$ 
 \ee[Y_{i_1}\cdots Y_{i_n} Y_{j_1}\cdots Y_{j_n}] = 
 {\bf 1}_{\{ i_1=j_1, \ldots , i_n=j_n\}} 
. 
$$ 
 Hence 
\begin{eqnarray*}
\lefteqn{
\ee[J_n(f_n)J_n(g_n)]
}
\\
 & = &
 (n!)^2 
 \ee\left[
\sum_{0\leq i_1 < \cdots < i_n}
 f_n(i_1,\ldots ,i_n) Y_{i_1}\cdots Y_{i_n}
\sum_{0\leq j_1 < \cdots < j_n} 
 g_n(j_1,\ldots ,j_n) Y_{j_1}\cdots Y_{j_n}
\right]
\\
 & = &
 (n!)^2 
 \sum_{0 \leq i_1 < \cdots < i_n, \ 
 0\leq j_1 < \cdots < j_n} 
 f_n(i_1,\ldots ,i_n)
 g_n(j_1,\ldots ,j_n)
 \ee[Y_{i_1}\cdots Y_{i_n} Y_{j_1}\cdots Y_{j_n}] 
\\
 & = &
 (n!)^2 \sum_{0 \leq i_1 < \cdots < i_n}
 f_n(i_1,\ldots ,i_n)
 g_n(i_1,\ldots ,i_n)
\\
 & = &
 n! \sum_{(i_1,\ldots ,i_n)\in \Delta_n
}
 f_n(i_1,\ldots ,i_n)
 g_n(i_1,\ldots ,i_n)
\\
 & = &
 n! \langle {\bf 1}_{\Delta_n} f_n , g_m \rangle_{\ell^2(\inte)^{\otimes n}}
.
\end{eqnarray*} 
 When $n < m$ and $(i_1,\ldots, i_n)\in \Delta_n$ and 
 $(j_1,\ldots, j_m)\in\Delta_m$ are two 
 sets of indices, there necessarily exists $k\in \{1,\ldots, m \}$ 
 such that $j_k\notin \{ i_1,\ldots, i_n \}$, hence 
$$\ee[Y_{i_1}\cdots Y_{i_n} Y_{j_1}\cdots Y_{j_m}] = 0 
, 
$$ 
 and this implies the orthogonality of $J_n (f_n)$ and $J_m (g_m)$.  
 The recurrence relation \eqref{recrel2} is a direct 
 consequence of \eqref{eg1}. 
 The isometry property \eqref{isois} of $J_n$ also follows 
 by induction from \eqref{iiiso} and the recurrence relation. 
\end{Proof} 
 If $f_n \in \ell^2(\inte^n)$ is not symmetric
 we let $J_n(f_n)=J_n(\tilde{f}_n)$, where
 $\tilde{f}_n$ is the symmetrization of $f_n$, defined as
$$\tilde{f}_n (i_1,\ldots ,i_n) 
 = \frac{1}{n!}
 \sum_{\sigma \in \Sigma_n} f(i_{\sigma (1)}, \ldots ,i_{\sigma_n}), 
 \qquad 
 i_1,\ldots ,i_n \in \inte^n, 
$$
 and $\Sigma_n$ is the set of all permutations of $\{1,\ldots ,n\}$.
 In particular, 
 if $(k_1, \ldots , k_n )\in \Delta_n$, 
 the symmetrization $\tilde{{\bf 1}}_{\{(k_1,\ldots ,k_n)\}}$ 
 of ${\bf 1}_{\{ (k_1,\ldots ,k_n) \}}$ in $n$ variables is given by 
$$ 
 {\bf \tilde{1}}_{\{(k_1,\ldots ,k_n)\}} 
 (i_1,\ldots ,i_n) 
 = \frac{1}{n!} {\bf 1}_{\{
 \{i_1,\ldots ,i_n\} 
 = \{k_1,\ldots ,k_n\}\}}, \quad 
 i_1,\ldots , i_n \in \inte, 
$$ 
 and 
$$
 J_n( 
 {\bf \tilde{1}}_{\{(k_1,\ldots ,k_n)\}})
 = Y_{k_1}\cdots Y_{k_n}. 
$$ 
\begin{lemma}
\label{lea2} 
 For all $n\geq 1$ we have 
$$\ee[J_n(f_n) \mid {\cal F}_k] = J_n (f_n {\bf 1}_{[0,k]^n}), 
 \qquad k\in\inte, \quad 
 f_n\in \ell^2 (\inte )^{\circ n}. 
$$ 
\end{lemma}
\begin{Proof} 
 This lemma can be proved in two ways, either as a consequence of 
 Proposition~\ref{isomp15} and Proposition~\ref{prop1.18} or via the 
 following direct argument, noting that for all 
 $m=0,\ldots ,n$ and $g_m \in \ell^2(\inte )^{\circ m}$ we have: 
\begin{eqnarray*} 
 \ee[(J_n(f_n) -J_n (f_n {\bf 1}_{[0,k]^n})) J_m (g_m {\bf 1}_{[0,k]^m}) ]
 & = & 
 {\bf 1}_{\{n=m\}} 
 n!\langle
 f_n(1-{\bf 1}_{[0,k]^n}) , g_m {\bf 1}_{[0,k]^m}\rangle_{\ell^2 (\inte^n)} 
\\ 
 & = & 
0 
, 
\end{eqnarray*} 
 hence $J_n (f_n {\bf 1}_{[0,k]^n})\in L^2(\Omega , {\cal F}_k)$, 
 and $J_n(f_n) -J_n (f_n {\bf 1}_{[0,k]^n})$ is orthogonal to 
 $L^2(\Omega , {\cal F}_k)$. 
\end{Proof} 
 In other terms we have 
$$\ee[ J_n(f_n)]=0, \qquad 
 f_n \in \ell^2(\inte)^{\circ n}, 
 \qquad 
 n\geq 1, 
$$ 
 the process $(J_n(f_n {\bf 1}_{[0,k]^n}))_{k\in\inte}$ is a discrete-time 
 martingale, and $J_n(f_n)$ is ${\cal F}_k$-measurable 
 if and only if $f_n {\bf 1}_{[0,k]^n}=f_n$, $0\leq k \leq n$. 
\section{Discrete structure equations}\index{structure equation!discrete time} 
\label{s3.1} 
 Assume now that the sequence $(Y_n)_{n\in\inte}$ 
 satisfies the discrete structure equation: 
\begin{equation}
\label{dse} 
Y_n^{2}=1+\varphi_n Y_n, \qquad n\in \inte,
\end{equation}
 where $(\varphi_n)_{n\in\inte}$ is an ${\cal F}_n$-predictable 
 process. 
 Condition \eqref{dse0} implies that 
$$\ee[Y_n^2 \mid {\cal F}_{n-1} ]=1, \qquad 
 n\in \inte 
, 
$$ 
 hence the hypotheses of the preceding sections are satisfied. 
 Since \eqref{dse} is a second order equation, there exists an 
 ${\cal F}_n$-adapted process $(X_n)_{n\in \inte}$ 
 of Bernoulli $\{-1,1\}$-valued random variables such that
\begin{equation} 
\label{yk} 
Y_n=\frac{\varphi_n}{2} + X_n \sqrt{1+\left(
 \frac{\varphi _n}{2}\right)^2},\qquad n\in\inte.
\end{equation} 
 Consider the conditional probabilities 
\begin{equation} 
\label{*p13} 
 p_n = \P(X_n=1 \mid {\cal F}_{n-1} ) \quad 
\mbox{and} \quad 
 q_n = \P(X_n=-1 \mid {\cal F}_{n-1} ), 
 \qquad n\in \inte.
\end{equation} 
 From the relation $\ee[Y_n \mid {\cal F}_{n-1} ]=0$, rewritten as 
$$p_n 
 \left( 
 \frac{\varphi_n}{2} + \sqrt{1+\left(
 \frac{\varphi _n}{2}\right)^2}
 \right) 
 + 
 q_n 
 \left( 
 \frac{\varphi_n}{2} - \sqrt{1+\left(
 \frac{\varphi _n}{2}\right)^2} 
 \right) 
 = 0, 
 \qquad 
 n\in \inte, 
$$ 
 we get 
\begin{equation} 
\label{yk01} 
 p_n = \frac{1}{2}
 \left( 
 1 -{\frac{\varphi _n}{\sqrt{ 4 + \varphi_n^{2}}}} 
 \right) 
,
 \qquad 
 q_n = {\frac{1}{2}} \left( 
 1 +{\frac{\varphi_n}{\sqrt{ 4 + \varphi_n^{2} }}} 
 \right) 
, 
\end{equation} 
 and 
$$ 
 \varphi_n = \sqrt{\frac{q_n}{p_n}} - \sqrt{\frac{p_n}{q_n}} 
 = \frac{q_n-p_n}{\sqrt{p_nq_n}}, 
 \qquad n\in \inte, 
$$ 
 hence 
$$ 
 Y_n 
 = {\bf 1}_{\{X_n=1\}} \sqrt{\frac{q_n}{p_n}} 
 - {\bf 1}_{\{X_n=-1\}} \sqrt{\frac{p_n}{q_n}}, 
 \qquad n\in \inte 
. 
$$ 
 Letting 
$$Z_n = \frac{X_n+1}{2} \in \{ 0,1\}, \qquad n\in\inte 
, 
$$ 
 we also have the relations 
\begin{equation}
\label{ch1}
Y_n = \frac{q_n-p_n+X_n}{2\sqrt{p_nq_n}} 
 = \frac{Z_n-p_n}{\sqrt{p_nq_n}}, 
 \quad n\in\inte 
, 
\end{equation} 
 which yield 
$${\cal F}_n = \sigma (X_0,\ldots ,X_n) 
 = \sigma (Z_0,\ldots ,Z_n)
, \qquad 
 n\in\inte 
. 
$$ 
\begin{remark} 
\label{rk01} 
 In particular, one can take $\Omega = \{-1,1\}^\inte$ and 
 construct the Bernoulli process $(X_n)_{n\in\inte}$ 
 as the sequence of canonical projections on 
 $\Omega = \{-1,1\}^\inte$ 
 under a countable product $\P$ of Bernoulli measures on $\{-1,1\}$. 
 In this case the sequence $(X_n)_{n\in\inte}$ can be viewed as the 
 dyadic expansion of $X(\omega ) \in [0,1]$ defined as: 
$$ 
X(\omega ) = \sum_{n=0}^\infty 
 \frac{1}{2^{n+1}} 
 X_n (\omega ) 
. 
$$ 
 In the symmetric case $p_k=q_k=1/2$, $k\in \inte$,
 the image measure of $\P$ by the mapping $\omega \mapsto X(\omega )$ 
 is the Lebesgue measure on $[0,1]$, 
 see \cite{MR1764269} for the non-symmetric case. 
\end{remark} 
\section{Chaos representation}\index{chaos representation!discrete time} 
\label{s3.2} 
 From now on we assume that the sequence $(p_k)_{k\in\inte}$ defined in 
 \eqref{*p13} is deterministic, which implies that the random variables 
 $(X_n)_{n\in\inte}$ are independent. 
 Precisely, $X_n$ will be  
 constructed as the canonical projection 
 $X_n:\Omega \to \{-1,1\}$ on $\Omega = \{-1,1\}^{\inte}$ 
 under the measure $\P$ given on cylinder sets by 
$$\P(\{\epsilon_0, \ldots, \epsilon_n \} \times \{ -1 , 1 \}^{\inte}) 
 = \prod_{k=0}^n 
 p_k^{(1+\varepsilon_k)/2} 
 q_k^{(1-\varepsilon_k)/2} 
, 
 \qquad \{ \epsilon_0, \ldots, \epsilon_n \} \in \{-1,1\}^{n+1} 
. 
$$ 
 The sequence $(Y_k)_{k\in\inte}$ can be constructed as 
 a family of independent random variables given by 
$$ 
 Y_n 
 = 
 \frac{\varphi_n}{2} + X_n \sqrt{1+\left(
 \frac{\varphi _n}{2}\right)^2},\qquad n\in\inte, 
$$ 
 where the sequence $(\varphi_n)_{n\in\inte}$ is deterministic. 
 In this case, all spaces 
 $L^r (\Omega, {\cal F}_n)$, $r\geq 1$, 
 have finite dimension $2^{n+1}$, with basis
\begin{eqnarray*} 
\lefteqn{ 
\left\{ {\bf 1}_{\{ 
 Y_0 = \epsilon_0, \ldots, 
 Y_n = \epsilon_n 
 \}}
 \ : \ (\epsilon_0,\ldots, \epsilon_n) \in 
 \prod_{k=0}^n 
 \left\{
 \sqrt{\frac{q_k}{p_k}}
 , -\sqrt{\frac{p_k}{q_k}} \right\} \right\} 
} 
\\ 
 & = & 
\left\{ {\bf 1}_{\{ 
 X_0 = \epsilon_0, \ldots, 
 X_n = \epsilon_n 
 \}}
 \ : \ (\epsilon_0,\ldots, \epsilon_n) \in 
 \prod_{k=0}^n 
 \left\{
 -1 , 1 
  \right\} \right\} 
. 
\end{eqnarray*} 
 An orthogonal basis of $L^r (\Omega, {\cal F}_n)$ is given by 
$$\left\{ Y_{k_1}\cdots Y_{k_l}
 = J_l ( 
 {\bf \tilde{1}}_{\{(k_1,\ldots ,k_l)\}}) 
  \ : \ 0\leq k_1 < \cdots < k_l \leq n,
 \ l=0,\ldots , n+1
 \right\} 
. 
$$
 Let 
\begin{equation} 
\label{defsn} 
 S_n = \sum_{k=0}^n \frac{1+X_k}{2} 
 = \sum_{k=0}^n Z_k, \quad n\in\inte, 
\end{equation} 
 denote the random walk associated to $(X_k)_{k\in \inte}$.
 If $p_k=p$, $k\in \inte$, then 
\begin{equation} 
\label{krw} 
J_n( {\bf 1}_{[0,N]}^{\circ n}) = K_n(S_N;N+1,p) 
\end{equation} 
 coincides with the Krawtchouk polynomial $K_n (\cdot ;N+1,p)$ of order $n$ 
 and parameter $(N+1,p)$, evaluated at $S_N$, cf. \cite{prisch2}.
\bigskip

\noindent 
 Let now ${\cal H}_0=\real$ and let
 ${\cal H}_n$ denote the subspace of $L^2(\Omega )$
 made of integrals of order $n\geq 1$, and called chaos 
 of order $n$: 
$${\cal H}_n = \{ J_n (f_n) \ : \
 f_n \in \ell^2(\inte)^{\circ n}\} 
.
$$ 
 The space of ${\cal F}_n$-measurable random variables is denoted by $L^0 (\Omega , {\cal F}_n)$. 
\begin{lemma} 
\label{ll1.1} 
 For all $n\in \inte$ we have 
\begin{equation} 
\label{ch}
 L^0 (\Omega , {\cal F}_n) \subset {\cal H}_0\oplus \cdots \oplus {\cal H}_{n+1}. 
\end{equation}
\end{lemma} 
\begin{Proof} 
 It suffices to note that ${\cal H}_l \cap L^0 (\Omega ,{\cal F}_n)$ 
 has dimension ${n+1\choose l}$, 
 $1 \leq l\leq n+1$. 
 More precisely it is generated by the orthonormal basis 
$$\left\{ Y_{k_1}\cdots Y_{k_l} = 
 J_l ( 
 {\bf \tilde{1}}_{\{(k_1,\ldots ,k_l)\}}) 
 \ : \ 0\leq k_1 < \cdots < k_l \leq n 
 \right\} 
, 
$$ 
 since any element $F$ of  ${\cal H}_l \cap L^0 (\Omega ,{\cal F}_n)$ 
 can be written as $F = J_l (f_l {\bf 1}_{[0,n]^l} )$, hence 
$$ 
 L^0 (\Omega , {\cal F}_n) = ( 
 {\cal H}_0\oplus \cdots \oplus {\cal H}_{n+1} ) 
 \bigcap L^0 (\Omega , {\cal F}_n) 
. 
$$ 
\end{Proof} 
\noindent 
 Alternatively, Lemma~\ref{ll1.1} can be proved by noting that 
$$J_n (f_n {\bf 1}_{[0,N]^n}) = 0, 
 \qquad 
 n>N+1, 
 \quad 
 f_n \in \ell^2(\inte)^{\circ n} 
, 
$$ 
 and as a consequence, any $F\in L^0 (\Omega , {\cal F}_N)$ 
 can be expressed as 
$$ 
 F = \ee[F] 
 + 
 \sum_{n=1}^{N+1} 
 J_n (f_n {\bf 1}_{[0,N]^n}) 
. 
$$ 
\begin{definition} 
 Let ${\cal S}$ denote the linear space 
 spanned by multiple stochastic integrals, i.e. 
\begin{equation} 
\label{defs} 
{\cal S} = \Vect \left\{ 
 \bigcup_{n=0}^\infty {\cal H}_n 
 \right\} 
 = \left\{
 \sum_{k=0}^n 
 J_k (f_k) \ : \
 f_k\in \ell^2(\inte )^{\circ k}
 , \ k=0, \ldots , n,\ n\in\inte\right\}. 
\end{equation} 
\end{definition} 
 The completion of ${\cal S}$ in $L^2(\Omega )$ 
 is denoted by the direct sum
$$\bigoplus_{n=0}^\infty
 {\cal H}_n 
. 
$$ 
 The next result is the chaos representation property for 
 Bernoulli processes, which is analogous to the Walsh decomposition, 
 cf. \cite{leitz}. 
 This property is obtained under the assumption that the sequence 
 $(X_n)_{n\in\inte}$ is i.i.d. 
 See \cite{emerycrp} for other instances of the chaos representation 
 property without this independence assumption. 
\begin{prop} 
\label{chaos}
 We have the identity
$$L^2(\Omega ) = \bigoplus_{n=0}^\infty
 {\cal H}_n.$$
\end{prop}
\begin{Proof} 
 It suffices to show that ${\cal S}$ is dense in $L^2(\Omega )$.
 Let $F$ be a bounded random variable. 
 Relation \eqref{ch} of Lemma~\ref{ll1.1} shows 
 that $\ee[F \mid {\cal F}_n] \in {\cal S}$. 
 The martingale convergence theorem, 
 cf. e.g. Theorem~27.1 in \cite{jacodprotterbk}, 
 implies that $(\ee[F \mid {\cal F}_n])_{n\in\inte}$ converges
 to $F$ a.s., hence every bounded
 $F$ is the $L^2 (\Omega )$-limit of a sequence in ${\cal S}$.
 If $F\in L^2(\Omega )$ is not bounded, $F$ is the limit
 in $L^2(\Omega )$ of the sequence
 $( {\bf 1}_{\{\vert F\vert\leq n\}} F)_{n\in\inte}$
 of bounded random variables.
\end{Proof} 
\noindent 
 As a consequence of Proposition~\ref{chaos},
 any $F\in L^{2}(\Omega ,\P)$ has a unique decomposition
\[
F= \ee[F] 
 + \sum_{n=1}^{\infty }J_{n}(f_{n}),\qquad f_{n}\in l^{2}(\inte)^{\circ n},\ n\in \inte,
\] 
 as a series of multiple stochastic integrals. 
 Note also that the statement of Lemma~\ref{ll1.1} is sufficient 
 for the chaos representation property to hold. 
\section{Gradient Operator}\index{finite difference gradient!discrete time}
\label{s4} 
\begin{definition}
 We densely define the linear gradient operator 
$$D: {\cal S} \longrightarrow L^{2}(\Omega \times \inte) 
$$ 
 by 
\[
D_{k}J_n(f_n)
 = n J_{n-1} (f_n(*,k) {\bf 1}_{\Delta_n}(*,k)), 
\]
 $k\in \inte$, 
 $f_n\in \ell^2(\inte)^{\circ n}$ 
 $n\in\inte$. 
\end{definition} 
 Note that for all 
 $k_1,\ldots ,k_{n-1} , k \in \inte$ 
 we have 
$$ 
{\bf 1}_{\Delta_n} ( k_1,\ldots ,k_{n-1} , k) 
 = 
{\bf 1}_{\{ k \notin ( k_1,\ldots , k_{n-1} ) \}} 
{\bf 1}_{\Delta_{n-1}} ( k_1,\ldots ,k_{n-1}) 
, 
$$ 
 hence we can write 
\[
 D_{k}J_n(f_n)
 = n J_{n-1} (f_n(*,k) {\bf 1}_{\{ k \notin * \}}), \quad
 k\in \inte, 
\]
 where in the above relation, ``$*$'' denotes the 
 first $k-1$ variables 
 $( k_1,\ldots , k_{n-1} )$ of $f_n( k_1,\ldots , k_{n-1} , k)$. 
 We also have 
 $D_k F = 0$ whenever $F\in {\cal S}$ is ${\cal F}_{k-1}$-measurable. 
\\ 

 On the other hand, $D_k$ is a continuous operator on the chaos 
 ${\cal H}_n$ since 
\begin{eqnarray} 
\nonumber 
 \Vert D_k J_n (f_n) \Vert_{L^2(\Omega )}^2 
 & = & 
 n^2 
 \Vert J_{n-1} (f_n(*,k)) \Vert_{L^2(\Omega )}^2 
\\ 
\label{frm11} 
 & = & 
 n n! 
 \Vert f_n(*,k) \Vert_{\ell^2(\inte^{\otimes (n-1)} )}^2 
, 
 \qquad 
 f_n \in \ell^2(\inte^{\otimes n} ), 
 \quad 
 k\in\inte. 
\end{eqnarray} 
 The following result gives the probabilistic interpretation of $D_k$
 as a finite difference operator. Given 
$$\omega = (\omega_0,\omega_1 , \ldots ) \in \{-1,1\}^\inte, 
$$ 
 let 
$$\omega_+^k = (\omega_0,\omega_1,\ldots ,\omega_{k-1}, +1,\omega_{k+1}, \ldots )$$ 
 and 
$$\omega_-^k = (\omega_0,\omega_1,\ldots ,\omega_{k-1}, -1,\omega_{k+1}, \ldots ) 
. 
$$ 
\begin{prop} 
\label{fimdtl} 
 We have for any $F\in {\cal S}$: 
\begin{equation} 
\label{natdef} 
 D_k F (\omega ) = \sqrt{p_kq_k}
 (F (\omega_+^k) - F (\omega_-^k ) ) , \quad k\in \inte 
. 
\end{equation} 
\end{prop}
\begin{Proof} 
 We start by proving the above statement for an ${\cal F}_n$-measurable 
 $F\in {\cal S}$. 
 Since $L^0(\Omega , {\cal F}_n)$ is finite dimensional it suffices to consider 
$$F= Y_{k_1}\cdots Y_{k_l} = f(X_0,\ldots ,X_{k_l}), 
$$ 
 with from \eqref{ch1}:
$$f(x_0,\ldots ,x_{k_l})
 = \frac{1}{2^l}
 \prod_{i=1}^{l}
 \frac{q_{k_i}-p_{k_i}+x_{k_i}}{\sqrt{p_{k_i}q_{k_i}}}.
$$ 
 First we note that from \eqref{ch} we have for 
 $(k_1,\ldots , k_n)\in \Delta_n$: 
\begin{eqnarray} 
\nonumber 
 D_{k}\left( Y_{k_1}\cdots Y_{k_n} \right)
 & = & 
 D_{k} J_n( 
 {\bf \tilde{1}}_{\{(k_1,\ldots ,k_n)\}})
\\ 
\nonumber 
 & = &  
 n J_{n-1} ( 
 {\bf \tilde{1}}_{\{(k_1,\ldots ,k_n)\}}(*,k))
\\ 
\nonumber 
 & = &  
 \frac{1}{(n-1)!} 
 \sum_{i=1}^{n} 
 {\bf 1}_{\{k_i\}}(k) 
 \sum_{(i_1,\ldots , i_{n-1})\in \Delta_{n-1}} 
 {\bf \tilde{1}}_{\{\{i_1,\ldots ,i_{n-1}\} = 
 \{ k_1,\ldots ,k_{i-1},k_{i+1}, \ldots , k_n 
 \} \} 
 }
\\ 
\nonumber 
 & = & \sum_{i=1}^{n} 
 {\bf 1}_{\{k_i\}}(k) 
 J_{n-1} ( 
 {\bf \tilde{1}}_{\{(k_1,\ldots ,k_{i-1},k_{i+1}, \ldots , k_n)\}} )
\\ 
\label{opl} 
 & = & 
 {\bf 1}_{\{k_1,\ldots ,k_n\}}(k) 
 \prod_{i=1 \atop k_i\not= k}^{n} Y_{k_i}.
\end{eqnarray} 
 If $k \notin \{k_1,\ldots ,k_l \}$ 
 we clearly have $F(\omega_+^k ) =F(\omega_-^k)=F(\omega )$, 
 hence 
$$\sqrt{p_kq_k}
 (F (\omega_+^k ) - F (\omega_-^k ) ) =0 
 = 
 D_kF (\omega ) 
. 
$$ 
 On the other hand if $ k \in \{k_1,\ldots ,k_l \}$ 
 we have 
$$F (\omega_+^k ) = 
 \sqrt{\frac{q_k}{p_k}} 
 \prod_{i=1 \atop k_i\not= k}^{l}
 \frac{q_{k_i}-p_{k_i}+\omega_{k_i}}{2 \sqrt{p_{k_i}q_{k_i}}}, 
$$ 
$$ 
 F (\omega_-^k ) 
 = 
 - \sqrt{\frac{p_k}{q_k}} 
 \prod_{i=1 \atop k_i\not= k}^{l}
 \frac{q_{k_i}-p_{k_i}+\omega_{k_i}}{2\sqrt{p_{k_i}q_{k_i}}}, 
$$ 
 hence from \eqref{opl} we get 
\begin{eqnarray*} 
\sqrt{p_kq_k}(F(\omega_+^k ) -F (\omega_-^k ) ) 
& = & 
 \frac{1}{2^{l-1}} 
 \prod_{i=1 \atop k_i\not= k}^{l}
 \frac{q_{k_i}-p_{k_i}+\omega_{k_i}}{\sqrt{p_{k_i}q_{k_i}}} 
\\ 
& = & 
 \prod_{i=1 \atop k_i\not= k}^{l} Y_{k_i} (\omega ) 
\\ 
 & = & 
 D_{k}\left( Y_{k_1}\cdots Y_{k_l} \right) (\omega ) 
\\ 
 & = & D_{k} F (\omega ) 
. 
\end{eqnarray*} 
 In the general case, $J_l(f_l )$ is the $L^2$-limit 
 of the sequence $\ee[J_l(f_l ) \mid {\cal F}_n] = J_l ( f_l {\bf 1}_{[0,n]^l})$ 
 as $n$ goes to infinity, and since from 
 \eqref{frm11} the operator $D_k$ is continuous 
 on all chaoses ${\cal H}_n$, $n\geq 1$, we have 
\begin{eqnarray*} 
 D_k F
 & = & 
 \lim_{n\to \infty} D_k \ee[F \mid {\cal F}_n] 
 \\ 
 & = & 
 \lim_{n\to \infty} 
 (\ee[F \mid {\cal F}_n] (\omega_+^k) - \ee[F \mid {\cal F}_n] (\omega_-^k ) )
\\ 
 & = & 
 \sqrt{p_kq_k} 
 (F (\omega_+^k) - F (\omega_-^k ) )
, \quad k\in \inte 
. 
\end{eqnarray*} 
\end{Proof} 
 The next property follows immediately from Proposition~\ref{fimdtl} . 
\begin{corollary} 
\label{cormes} 
 A random variable $F:\Omega \to \real$ 
 is ${\cal F}_n$-measurable if and only if 
$$D_k F = 0$$ 
 for all $k>n$. 
\end{corollary} 
 If $F$ has the form $F = f(X_0,\ldots ,X_n)$, 
 we may also write 
$$ 
D_k F = \sqrt{p_kq_k}
 (F^+_k - F^-_k ) , \qquad k\in \inte, 
$$ 
 with 
$$F_k^+ = f(X_0,\ldots ,X_{k-1},+1,X_{k+1}, \ldots ,X_n), 
$$ 
 and 
$$
F_k^- = 
 f(X_0,\ldots ,X_{k-1},-1,X_{k+1}, \ldots ,X_n). 
$$ 
 The gradient $D$ can also be expressed as
$$ 
 D_k F (S_\cdot ) = \sqrt{p_kq_k}
 \left(
 F 
 \left( 
 S_\cdot + {\bf 1}_{\{X_k=-1\}} {\bf 1}_{\{k\leq \cdot \}} 
 \right) 
 - 
 F 
 \left( 
 S_\cdot - {\bf 1}_{\{X_k=1\}} {\bf 1}_{\{k\leq \cdot \}} 
 \right) 
 \right) 
 ,
$$
 where $F (S_\cdot )$ is an informal notation for 
 the random variable $F$ estimated on a given path 
 of $(S_n)_{n\in\inte}$ defined in \eqref{defsn} 
 and $S_\cdot + {\bf 1}_{\{X_k=\mp 1 \}} {\bf 1}_{\{k\leq \cdot \}}$ 
 denotes the path of $(S_n)_{n\in\inte}$ perturbed by 
 forcing $X_k$ to be equal to $\pm 1$. 
\\ 

 We will also use the gradient $\nabla_k$ defined as 
\begin{equation} 
\label{mod2} 
\nabla_k F =
 X_k
 \left(
 f(X_0,\ldots ,X_{k-1}, -1,X_{k+1}, \ldots ,X_n)
 - f(X_0,\ldots ,X_{k-1},1,X_{k+1}, \ldots ,X_n)
 \right),
\end{equation} 
 $k\in \inte$, 
 with the relation 
$$ 
 D_k = - X_k \sqrt{p_kq_k} \nabla_k 
, 
 \qquad 
 k\in \inte 
, 
$$ 
 hence $\nabla_kF$ coincides with 
 $D_kF$ after squaring and multiplication by $p_kq_k$. 
 From now on, $D_k$ denotes the finite difference 
 operator which is extended to any $F:\Omega \to \real$ 
 using Relation~\eqref{natdef}. 
 The $L^2$ domain of $D$ is naturally defined as 
 the space of functionals $F$ such that $\ee[\Vert DF\Vert^2_{\ell^2 (\inte)}] 
 < \infty$, or equivalently 
$$\sum_{n=0}^\infty n n! \Vert f_n \Vert_{\ell^2 (\inte^n)}^2 < \infty 
, 
$$ 
 if $F = \sum_{n=0}^\infty J_n ( f_n )$. 
 The following is the product rule for the operator $D$.
\begin{prop} 
\label{chnrle} 
 Let $F,G:\Omega \rightarrow \real$.
 We have
\begin{eqnarray*}
D_{k}(FG)& =& FD_{k}G+GD_{k}F-\frac{X_k}{\sqrt{p_kq_k}} D_{k}FD_{k}G,
 \qquad k\in \inte. 
\end{eqnarray*}
\end{prop}
\begin{Proof} 
 Let $F_+^k (\omega ) = F (\omega_+^k)$, $F_-^k (\omega ) = F (\omega_-^k)$, $k\geq 0$.
 We have 
\begin{eqnarray*}
D_k(FG) & = &  \sqrt{p_kq_k} (F_+^k G_+^k - F_-^k G_-^k )
\\
& = & {\bf 1}_{\{X_k=-1\}}
 \sqrt{p_kq_k}
 \left(
 F(G_+^k -G)+
 G(F_+^k -F)
 +(F_+^k -F)(G_+^k -G )\right)
\\
& & + {\bf 1}_{\{X_k=1\}} \sqrt{p_kq_k}
 \left( F(G-G_-^k )+
 G(F-F_-^k )
 -(F-F_-^k )(G-G_-^k )\right)
\\
& = &
 {\bf 1}_{\{X_k=-1\}}
\left(
FD_kG+GD_kF+\frac{1}{\sqrt{p_kq_k}}D_kFD_kG\right)
\\
& & + {\bf 1}_{\{X_k=1\}} \left(
 FD_kG+GD_kF-\frac{1}{\sqrt{p_kq_k}}D_kFD_kG\right).
\end{eqnarray*}
\end{Proof} 
\section{Clark Formula and Predictable Representation}\index{Clark formula!discrete time} \index{predictable representation!discrete time} 
\label{s5} 
 In this section we prove a 
 predictable representation formula for the functionals of $(S_n)_{n\geq 0}$ 
 defined in \eqref{defsn}. 
\begin{prop} 
\label{propclk} 
 For all $F \in {\cal S}$ we have 
\begin{eqnarray} 
\label{clk} 
F & = & \ee[F]+\sum_{k=0}^{\infty }\ee[D_{k}F \mid {\cal F}_{k-1}]Y_{k} 
\\ 
\nonumber 
 & = & 
 \ee[F]+\sum_{k=0}^{\infty } Y_{k} D_{k} \ee[F \mid {\cal F}_k ] 
. 
\end{eqnarray} 
\end{prop}
\begin{Proof}
 The formula is obviously true for $F = J_0(f_0)$. 
 Given $n\geq 1$, as a consequence of Proposition~\ref{prop1.18} above and
 Lemma~\ref{lea2} we have: 
\begin{eqnarray*} 
 J_n (f_n) 
& = & n\sum_{k=0}^{\infty }J_{n-1}(f_{n}( * ,k) {\bf 1}_{[0,k-1]^{n-1}}( *
))Y_{k}
\\
& = & n\sum_{k=0}^{\infty }J_{n-1}(f_{n}( * ,k) {\bf 1}_{\Delta_n}(*,k) {\bf 1}_{[0,k-1]^{n-1}}( *
))Y_{k}
\\
& = & n\sum_{k=0}^{\infty }\ee[J_{n-1}(f_{n}( * ,k) {\bf 1}_{\Delta_n}(*,k) ) \mid {\cal F}_{k-1}]
 Y_{k}
\\
& = & \sum_{k=0}^{\infty }\ee[D_kJ_n(f_{n} ) \mid {\cal F}_{k-1}]
 Y_{k}, 
\end{eqnarray*} 
 which yields \eqref{clk} for $F = J_n(f_n)$, since $\ee[ J_n(f_n)]=0$. 
 By linearity the formula is established for $F\in {\cal S}$. 
\end{Proof} 
 Although the operator 
 $D$ is unbounded we have the following result, 
 which states the boundedness of the operator that maps 
 a random variable to the unique process involved in its predictable 
 representation. 
\begin{lemma} 
\label{lbdd} 
 The operator 
\begin{align*} 
L^{2}(\Omega ) & \longrightarrow L^{2}(\Omega \times \inte) 
\\ 
 F & \mapsto (\ee[D_k F \mid {\cal F}_{k -1}])_{k\in\inte}
\end{align*} 
 is bounded with norm equal to one. 
\end{lemma} 
\begin{Proof} 
 Let $F\in {\cal S}$. 
 From Relation~\eqref{clk} and 
 the isometry formula \eqref{lkoi} for the stochastic 
 integral operator $J$ we get 
\begin{eqnarray} 
\label{bndd} 
\Vert \ee[D_{\cdot }F \mid {\cal F}_{\cdot -1}]\Vert _{L^{2}(\Omega \times \inte)}^{2} 
 & = & 
 \Vert F-\ee[F]\Vert_{L^{2}(\Omega )}^{2}  
\\  
\nonumber 
 & \leq &  
 \Vert F-\ee[F]\Vert _{L^{2}(\Omega )}^{2} + (\ee[F])^2 
\\ 
\nonumber 
 & = & \Vert F\Vert _{L^{2}(\Omega )}^{2}, 
\end{eqnarray} 
 with equality in case $F = J_1 (f_1)$. 
\end{Proof} 
 As a consequence of Lemma~\ref{lbdd} we have 
 the following corollary. 
\begin{corollary} 
 The Clark formula of 
 Proposition~\ref{propclk} 
 extends to any $F\in L^2 (\Omega )$. 
\end{corollary} 
\begin{Proof} 
 Since $F\mapsto \ee[D_{\cdot }F \mid {\cal F}_{\cdot -1}]$ is bounded 
 from Lemma~\ref{lbdd}, the Clark formula extends to $F\in L^{2}(\Omega )$ 
 by a standard Cauchy sequence argument.
 For the second identity we use the relation 
$$\ee[D_k F \mid {\cal F}_{k-1} ] = D_k \ee[F \mid {\cal F}_k ] 
$$ 
 which clearly holds since $D_k F$ is independent 
 of $X_k$, $k\in \inte$. 
\end{Proof} 
 Let us give a first elementary application 
 of the above construction to the proof of a Poincar\'e inequality 
 on Bernoulli space. We have 
\begin{eqnarray*} 
\var (F) & = & \ee[ | F-\ee[F] |^2] 
\\ 
 & = & \ee\left[ 
 \left( 
 \sum_{k=0}^\infty 
 \ee[D_k F \mid {\cal F}_{k-1} ] Y_k \right)^2 
 \right] 
\\ 
 & = & \ee\left[ 
 \sum_{k=0}^\infty 
 (\ee[D_k F \mid {\cal F}_{k-1} ])^2 
 \right] 
\\ 
 & \leq & \ee\left[ 
 \sum_{k=0}^\infty 
 \ee[ | D_k F |^2 \mid {\cal F}_{k-1} ] 
 \right] 
\\ 
 & = & \ee\left[ 
 \sum_{k=0}^\infty 
 | D_k F |^2 
 \right], 
\end{eqnarray*} 
 hence 
$$\var (F) \leq \Vert DF\Vert_{L^2(\Omega \times \inte )}^2. 
$$ 
 More generally the Clark formula implies the following. 
\begin{corollary} 
\label{lemmaa} 
 Let $a\in \inte$ and $F\in L^2 (\Omega )$. 
 We have 
\begin{equation} 
\label{clark2a} 
F= \ee[F \mid {\cal F}_a] 
 + \sum_{k=a+1}^\infty \ee[D_k F \mid {\cal F}_{k-1} ] Y_k, 
\end{equation} 
 and 
\begin{equation} 
\label{clark3} 
 \ee[F^2] =  \ee[(\ee[F \mid {\cal F}_a])^2] 
 + \ee\left[ \sum_{k=a+1}^\infty (\ee[D_k F \mid {\cal F}_{k-1}])^2 \right]. 
\end{equation} 
\end{corollary} 
\begin{Proof} 
 From Proposition~\ref{isomp15} and the Clark formula 
 \eqref{clk} of Proposition~\ref{propclk} we have 
$$ 
\ee[F \mid {\cal F}_a ] = 
 \ee[F]+\sum_{k=0}^a \ee[D_k F \mid {\cal F}_{k-1} ]Y_k 
, 
$$ 
 which implies (\ref{clark2a}). Relation (\ref{clark3}) is an immediate consequence of 
 (\ref{clark2a}) and the isometry property of $J$. 
\end{Proof} 
 As an application of the Clark formula of Corollary~\ref{lemmaa} 
 we obtain the following predictable representation property for 
 discrete-time martingales. 
\begin{prop} 
\label{martrepr} 
 Let $(M_n)_{n\in \inte}$ be a martingale in $L^2(\Omega)$ with respect 
 to $({\cal F}_n)_{n\in \inte}$. 
 There exists a predictable process $(u_k)_{k\in \inte}$ locally in $L^2(\Omega \times \inte)$, 
 (i.e. $u(\cdot ) {\bf 1}_{[0,N]} (\cdot ) \in L^2(\Omega \times \inte)$ 
 for all $N>0$) 
 such that 
\begin{equation} 
\label{ghl0}
 M_n = M_{-1} + \sum_{k=0}^n 
 u_k Y_k, \qquad n\in \inte.
\end{equation} 
\end{prop}
\begin{Proof} 
 Let $k\geq 1$. From 
 Corollaries~\ref{cormes} and \ref{lemmaa} 
 we have: 
\begin{eqnarray*} 
 M_k & = & 
 \ee[ M_k \mid {\cal F}_{k-1}] 
 + 
 \ee[D_k M_k \mid {\cal F}_{k-1} ] Y_k 
\\ 
& = & 
 M_{k-1} 
 + 
 \ee[D_k M_k \mid {\cal F}_{k-1} ] Y_k 
, 
\end{eqnarray*} 
 hence it suffices to let 
$$u_k = \ee[D_k M_k \mid {\cal F}_{k-1}], 
 \qquad k\geq 0,$$
 to obtain 
$$ 
 M_n = M_{-1} + \sum_{k=0}^n M_k-M_{k-1} 
 = M_{-1} + \sum_{k=0}^n u_k Y_k. 
$$ 
\end{Proof} 
\section{Divergence Operator}\index{divergence operator!discrete time} 
\label{s6} 
 The divergence operator $\delta$ is introduced as the adjoint of $D$. 
 Let ${\cal U}\subset L^2 (\Omega \times \inte )$ 
 be the space of processes defined as 
$${\cal U} = \left\{
 \sum_{k=0}^n 
 J_k (f_{k+1} (*,\cdot ) ),
 \quad
 f_{k+1} \in \ell^2(\inte )^{\circ k}\otimes \ell^2 (\inte), \
 \ k=, n \in\inte \right\} 
. 
$$ 
\begin{definition}
 Let $\delta : {\cal U} 
 \to L^2(\Omega )$ be the linear mapping defined on
${\cal U}$ 
 as
\[ 
 \delta ( u ) 
 = 
 \delta (J_{n}(f_{n+1}( * ,\cdot )))=J_{n+1}(\tilde{f}_{n+1}), \quad 
 f_{n+1}\in l^{2}(\inte)^{\circ n}\otimes l^{2}(\inte),
\]
 for $(u_k)_{k\in\inte}$ of the form 
$$ 
 u_k = 
 J_{n}(f_{n+1}( * , k )), 
 \qquad 
 k\in\inte 
, 
$$ 
 where $\tilde{f}_{n+1}$ denotes the symmetrization of $f_{n+1}$
 in $n+1$ variables, i.e. 
$$\tilde{f}_{n+1} (k_1,\ldots ,k_{n+1}) 
 = \frac{1}{n+1}
 \sum_{i=1}^{n+1}
 f_{n+1} (k_1,\ldots ,k_{k-1},k_{k+1}, \ldots ,k_{n+1},
 k_i ) 
. 
$$ 
\end{definition} 
 From Proposition~\ref{chaos}, ${\cal S}$ is dense in $L^2 (\Omega )$, 
 hence ${\cal U}$ is dense in $L^2 (\Omega \times \inte )$.  
\begin{prop} 
\label{prdual} 
 The operator $\delta$ is adjoint to $D$: 
$$\ee[\langle DF,u\rangle_{\ell^2 (\inte )} ]
 = \ee[F\delta (u) ], \quad F\in {\cal S}, \ u\in {\cal U}.
$$ 
\end{prop}
\begin{Proof}
 We consider $F= J_n(f_n)$ and $u_k = J_m(g_{m+1}(*,k))$, $k\in \inte$,
 where $f_n\in \ell^2(\inte)^{\circ n}$ and 
 $g_{m+1}\in \ell^2(\inte)^{\circ m}\otimes \ell^2(\inte )$.
 We have 
\begin{eqnarray*}
\lefteqn{ 
 \ee[\langle D_\cdot
 J_n(f_n), J_m (g_{m+1}(*,\cdot )) \rangle_{\ell^2 (\inte )}] 
 = 
 n \ee [ \langle 
 J_{n-1} (f_n ( * , \cdot ) ) 
 , 
 J_m (g_m ( * , \cdot ) ) 
 \rangle_{l^2 (\inte )} 
 ] 
 } 
\\
 & = &
 n {\bf 1}_{\{n-1=m \}}
\sum_{k=0}^\infty
\ee[J_{n-1}( f_n (*,k)  {\bf 1}_{\Delta_n }(*,k) )
 J_m(g_{m+1}(*,k) ) ] 
\\
 & = &
 n! {\bf 1}_{\{n-1=m \}}
\sum_{k=0}^\infty
 \langle
 {\bf 1}_{\Delta_n } (*,k) 
 f_n (*,k), g_{m+1} (*,k) 
 \rangle_{\ell^2 (\inte^{n-1})}
\\
 & = & 
 n! {\bf 1}_{\{n=m+1\}}
 \langle
 {\bf 1}_{\Delta_n }
 f_n, g_{m+1} \rangle_{\ell^2 (\inte^n)}
\\
 & = & 
 n! {\bf 1}_{\{n=m+1\}}
 \langle {\bf 1}_{\Delta_n}
 f_n, \tilde{g}_{m+1} \rangle_{\ell^2 (\inte^n)}
\\
 & = &
\ee[J_n(f_n)J_m( \tilde{g}_{m+1} ) ]
 \\ 
 & = & \ee [ \delta (u) F ]. 
\end{eqnarray*}
\end{Proof} 
 The next proposition shows in particular that $\delta$ coincides with 
 the stochastic integral operator $J$ on the square-summable 
 predictable processes. 
\begin{prop} 
\label{pkl} 
 The operator $\delta$ can be extended to $u\in L^2(\Omega \times \inte)$ 
 with 
\begin{equation} 
\label{dlta} 
\delta (u) = \sum_{k=0}^\infty u_k Y_k - \sum_{k=0}^\infty D_k u_k - \delta ( \varphi D u ), 
\end{equation} 
 provided that all series converges in $L^2 (\Omega )$, 
 where $(\varphi_k)_{k\in\inte}$ appears in the 
 structure equation \eqref{dse}.
 We also have for all $u \in {\cal U}$: 
\begin{equation} 
 \label{skois} 
 \ee[ | \delta (u)| ^2 ] = \ee[\Vert u \Vert_{\ell^2(\inte )}^2 ]
 + \ee\left[ \sum_{k,l=0\atop k \not= l}^\infty D_ku_l D_lu_k
  - \sum_{k=0}^\infty ( D_ku_k)^2 
   \right]. 
\end{equation} 
\end{prop} 
\begin{Proof} 
 Using the expression \eqref{eg1} of 
 $u_k = J_n (f_{n+1} (*,k))$ we have 
\begin{eqnarray*} 
\delta (u) 
 & = & 
 J_{n+1}(\tilde{f}_{n+1} ) 
\\ 
 & = & 
 \sum_{(i_1 , \ldots , i_{n+1} )\in \Delta_{n+1} } 
 \tilde{f}_{n+1} (i_1,\ldots ,i_{n+1} ) Y_{i_1}\cdots Y_{i_{n+1}} 
\\ 
 & = & 
 \sum_{k=0}^\infty 
 \sum_{(i_1 , \ldots , i_n )\in \Delta_n } 
 \tilde{f}_{n+1} (i_1,\ldots ,i_n , k ) Y_{i_1}\cdots Y_{i_n} Y_k 
\\ 
 & & 
 - 
 n 
 \sum_{k=0}^\infty 
 \sum_{(i_1 , \ldots , i_{n-1} )\in \Delta_{n-1} } 
 \tilde{f}_{n+1} (i_1,\ldots ,i_{n-1} , k , k ) Y_{i_1}\cdots Y_{i_{n-1}} | Y_k |^2 
\\ 
 & = & 
\sum_{k=0}^\infty 
 u_k Y_k 
 - 
 \sum_{k=0}^\infty D_k u_k | Y_k |^2 
\\ 
 & = & \sum_{k=0}^\infty 
 u_k Y_k 
 - 
 \sum_{k=0}^\infty D_k u_k 
 - 
 \sum_{k=0}^\infty \varphi_k D_k u_k Y_k. 
\end{eqnarray*} 
 Next, we note the commutation relation\footnote{See A.~Mantei, Masterarbeit ``Stochastisches Kalk\"ul in diskreter Zeit'', Satz 6.7, 2015.}  
\begin{eqnarray*} 
 D_k \delta (u) 
 & = & 
 D_k \left(
 \sum_{l=0}^\infty 
 u_l Y_l 
 - 
 \sum_{l=0}^\infty | Y_l |^2 D_l u_l
 \right)
 \\
 & = & 
 \sum_{l=0}^\infty 
 \left( Y_l D_k u_l + u_l D_k Y_l - \frac{X_k}{\sqrt{p_kq_k}} D_k u_l D_k Y_l
 \right) 
 \\
 & &
 - 
 \sum_{l=0}^\infty
 \left(
 | Y_l |^2 D_k D_l u_l
 + D_l u_l D_k | Y_l |^2 
 - \frac{X_k}{\sqrt{p_kq_k}} D_k | Y_l |^2 D_k D_l u_l
 \right)
 \\
 & = & 
 \delta ( D_k u)
 +
 u_k D_k Y_k - \frac{X_k}{\sqrt{p_kq_k}} D_k u_k D_k Y_k
 - D_k u_k D_k | Y_k |^2 
 \\
 & = & 
 \delta ( D_k u)
 + u_k -
 \left( \frac{X_k}{\sqrt{p_kq_k}} + 2 Y_k D_k Y_k -\frac{X_k}{\sqrt{p_kq_k}}
 D_k Y_k D_k Y_k 
 \right) D_k u_k 
 \\
 & = & 
 \delta ( D_k u) + u_k - 2 Y_k D_k u_k. 
\end{eqnarray*} 
 On the other hand, we have 
\begin{eqnarray*} 
 \delta ( {\bf 1}_{\{ k \} } D_k u_k ) 
 & = & 
 \sum_{l=0}^\infty 
 Y_l {\bf 1}_{\{ k \} } (l) D_k u_k 
 - 
 \sum_{l=0}^\infty | Y_l |^2 D_l ( {\bf 1}_{\{ k \} } (l) D_k u_k ) 
 \\
 & = & Y_k D_k u_k - | Y_k |^2 D_k D_k u_k 
 \\
 & = & Y_k D_k u_k, 
\end{eqnarray*} 
hence
\begin{eqnarray*} 
\Vert \delta (u) \Vert^2_{L^2(\Omega ) } 
& = & 
\ee[
  \langle u , D \delta (u) \rangle_{\ell^2(\inte )}
 ] 
\\
 & = & 
\ee \left[
  \sum_{k=0}^\infty
  u_k ( u_k + \delta ( D_k u ) - 2 Y_k D_k u_k ) 
 \right] 
\\
 & = & 
\ee[\Vert 
  u \Vert_{\ell^2(\inte )}^2 ]
 + \ee\left[ 
   \sum_{k,l=0}^\infty D_ku_l D_lu_k \right]
  - 2 \ee \left[ \sum_{k=0}^\infty u_k Y_k D_k u_k \right] 
\\
 & = & 
 \ee[\Vert 
 u \Vert_{\ell^2(\inte )}^2 ] 
  + \ee\left[ 
    \sum_{k,l=0}^\infty D_ku_l D_lu_k
    - 2 \sum_{k=0}^\infty ( D_ku_k )^2 
    \right], 
\end{eqnarray*} 
 where we used the equality
 \begin{eqnarray*}
  \ee \left[ u_k Y_k D_k u_k \right] 
   & = & 
  \ee \left[ p_k {\bf 1}_{\{ X_k = 1 \} } u_k (\omega_+^k)
    Y_k (\omega_+^k) D_k u_k 
  +
    q_k {\bf 1}_{\{ X_k = - 1 \} } u_k (\omega_-^k)
    Y_k (\omega_-^k) D_k u_k \right] 
   \\
   & = & 
   \sqrt{p_kq_k} \ee \left[
     ( {\bf 1}_{\{ X_k = 1 \} } u_k (\omega_+^k)
     - {\bf 1}_{\{ X_k = - 1 \} } u_k (\omega_-^k)
     )
     D_k u_k \right] 
    \\
    & = &
    \ee\left[ ( D_ku_k )^2 \right], \qquad k \in \inte. 
\end{eqnarray*} 
\end{Proof} 
 In the symmetric case $p_k=q_k=1/2$ we have $\varphi_k=0$, $k\in\inte$, and 
 $$ \delta (u) = \sum_{k=0}^\infty u_k Y_k - \sum_{k=0}^\infty D_k u_k. $$ 
 The last two terms in the right hand side of 
 \eqref{dlta} vanish when $(u_k)_{k\in\inte}$ is predictable, 
 and in this case the Skorohod isometry \eqref{skois} becomes
 the It\^o isometry as in the next proposition. 
\begin{corollary} 
 If $(u_k )_{k\in\inte}$ satisfies $D_ku_k = 0$, 
 i.e. $u_k$ does not depend on $X_k$, $k\in\inte$, then $\delta (u)$ 
 coincides with the (discrete time) stochastic integral 
\begin{equation} 
\label{iisoo} 
\delta (u)= \sum_{k=0}^{\infty }Y_{k}u_k 
,
\end{equation} 
 provided that the series converges in $L^2(\Omega )$. 
 If moreover $(u_k)_{k\in\inte}$ is predictable and square-summable 
 we have the isometry 
\begin{equation}
\label{iiso}
 \ee[\delta (u)^{2}]=\ee\left[\Vert u\Vert _{\ell^{2}(\inte )}^{2}\right], 
\end{equation} 
 and $\delta (u)$ coincides with $J (u)$ 
 on the space of predictable square-summable processes. 
\end{corollary} 
\section{Ornstein-Uhlenbeck Semi-Group and Process}\index{Ornstein-Uhlenbeck!semigroup}\index{Ornstein-Uhlenbeck!process}
\label{s7} 
 The Ornstein-Uhlenbeck operator $L$ is defined as $L = \delta D$, i.e. $L$ satisfies 
$$LJ_n ( f_n) = n J_n (f_n) , \qquad f_n \in \ell^2 (\inte)^{\circ n} 
. 
$$ 
\begin{prop} 
 For any $F\in {\cal S}$ we have 
$$LF = \delta D F = \sum_{k=0}^\infty Y_k (D_kF) 
 = \sum_{k=0}^\infty \sqrt{p_kq_k} Y_k (F_k^+ - F_k^- ) ,  
$$ 
\end{prop} 
\begin{Proof} 
 Note that $D_kD_kF=0$, $k\in\inte$, and 
 use Relation \eqref{dlta} of Proposition~\ref{pkl}. 
\end{Proof} 
 Note that $L$ can be expressed in other forms, for example 
$$ 
 LF = 
 \sum_{k=0}^\infty 
 \Delta_k F 
, 
$$ 
 where 
\begin{eqnarray*} 
 \Delta_k F 
 & = & 
 ( 
 {\bf 1}_{\{X_k =1\}} 
 q_k 
 (F (\omega ) - F (\omega_-^k ) ) 
 - 
 {\bf 1}_{\{X_k=-1\}} 
 p_k 
 (F (\omega_+^k) - F (\omega ) ) 
 ) 
\\ 
 & = & 
 F 
 - 
 ( 
 {\bf 1}_{\{X_k =1\}} 
 q_k 
 F (\omega_-^k ) 
 + 
 {\bf 1}_{\{X_k=-1\}} 
 p_k 
 F (\omega_+^k) 
 )
\\ 
 & = & 
 F 
 - 
 \ee 
 [ 
 F 
 \mid 
 {\cal F}_k^c 
 ] 
, \quad k\in \inte 
, 
\end{eqnarray*} 
 and 
 ${\cal F}_k^c$ is the $\sigma$-algebra generated by 
$$ 
 \{ X_l \ : \ l\not=k, \ l \in \inte \}. 
$$ 
 Let now $(P_t)_{t\in \real_+} = (e^{tL})_{t\in \real_+}$ 
 denote the semi-group associated to $L$ and defined as
$$P_t F = \sum_{n=0}^\infty e^{-nt} J_n(f_n),
 \qquad t\in \real_+,$$
 on 
 $\displaystyle F = \sum_{n=0}^\infty J_n(f_n) \in L^2(\Omega )$. 
 The next result shows that $(P_t)_{t\in \real_+}$
 admits an integral representation by a probability kernel.
 Let $q^N_t : \Omega \times \Omega \to \real_+$ be defined by 
$$ 
 q^N_t(\tilde{\omega} , \omega ) =
 \prod_{i=0}^{N}
 (1+e^{-t} Y_i( \omega ) Y_i( \tilde{\omega} ) ),
 \qquad \omega, \tilde{\omega} \in \Omega, 
 \quad t\in \real_+. 
$$ 
\begin{lemma}
\label{p11} 
 Let the probability kernel $Q_t (\tilde{\omega} , d\omega )$ be defined by 
$$\ee\left[ 
 \frac{dQ_t (\tilde{\omega} , \cdot )}{d\P} \Big| 
 {\cal F}_N \right] (\omega ) 
 = 
 q^N_t(\tilde{\omega} , \omega ) 
, \qquad N\geq 1, 
 \quad t\in \real_+. 
$$ 
 For $F\in L^2(\Omega, {\cal F}_N)$ we have
\begin{equation}
\label{repr}
 P_t F ( \tilde{\omega} ) =
 \int_\Omega F( \omega ) 
 Q_t ( \tilde{\omega} , d \omega ), 
 \qquad \tilde{\omega} \in \Omega, \quad n\geq N. 
\end{equation} 
\end{lemma} 
\begin{Proof}
 Since $L^2(\Omega ,{\cal F}_N)$ has finite
 dimension $2^{N+1}$,
 it suffices to consider functionals of the form
 $F = Y_{k_1}\cdots Y_{k_n}$ with $0\leq k_1<\cdots < k_n\leq N$.
 We have for $\omega \in \Omega$, $k\in \inte$:
\begin{eqnarray*}
\lefteqn{
\ee\left[Y_k(\cdot) (1+e^{-t} Y_k(\cdot ) Y_k(\omega ) )\right]
}
\\
& = &
 p_k
 \sqrt{\frac{q_k}{p_k}}
 \left(1+e^{-t} \sqrt{\frac{q_k}{p_k}} Y_k (\omega )\right)
 - q_k
 \sqrt{\frac{p_k}{q_k}}
 \left(1-e^{-t} \sqrt{\frac{p_k}{q_k}} Y_k (\omega )\right)
\\
 & = &
 e^{-t} Y_k(\omega ),
\end{eqnarray*} 
 which implies, by independence of the sequence $(X_k)_{k\in\inte}$, 
\begin{eqnarray*} 
\ee[Y_{k_1}\cdots Y_{k_n} q^N_t(\omega , \cdot )]
& = &
\ee\left[
 Y_{k_1}\cdots Y_{k_n}
 \prod_{i=1}^{N}
 (1+e^{-t} Y_{k_i}(\omega ) Y_{k_i} (\cdot ) )
 \right] 
\\
& = &
 \prod_{i=1}^{N}
\ee\left[
 Y_{k_i} (\cdot )
 (1+e^{-t} Y_{k_i} (\omega ) Y_{k_i} (\cdot ) )
 \right]
\\
& = &
 e^{-nt} Y_{k_1}(\omega )\cdots Y_{k_n}(\omega )
\\
& = & 
 e^{-nt} J_n( {\bf \tilde{1}}_{\{(k_1, \ldots , k_n )\}} ) (\omega )
\\
& = & 
 P_t J_n( 
 {\bf \tilde{1}}_{\{(k_1, \ldots , k_n )\}} ) (\omega )
\\ 
& = & P_t(Y_{k_1}\cdots Y_{k_n}) (\omega ). 
\end{eqnarray*}
\end{Proof} 
 Consider the $\Omega$-valued stationary process $(X(t))_{t\in \real_+}
 = ((X_k(t) )_{k\in\inte})_{t\in \real_+}$ 
 with independent components and distribution given by 
\begin{align} 
\label{a1.x} 
& \P(X_k (t) = 1 \mid X_k (0) = 1) = p_k+e^{-t}q_k,
\\ 
\nonumber 
\\ 
\label{a2.x} 
& \P(X_k (t) = -1 \mid X_k (0) = 1) = q_k-e^{-t}q_k,
\\ 
\nonumber 
\\ 
\label{a3.x} 
& \P(X_k (t) = 1 \mid X_k (0) = -1) = p_k-e^{-t}p_k, 
\\ 
\nonumber 
\\ 
\label{a4.x} 
& \P (X_k (t) = -1 \mid X_k (0) = -1) = q_k+e^{-t}p_k, 
\end{align} 
 $k\in\inte$, $t\in \real_+$. 
\begin{prop} 
 The process $(X (t))_{t\in \real_+} = ((X_k (t))_{k\in\inte})_{t\in \real_+}$ 
 is the Ornstein-Uhlenbeck process associated to $(P_t)_{t\in \real_+}$, 
 i.e. we have 
\begin{equation} 
\label{istheou} 
 P_t F = \ee[F(X (t) ) \mid X (0) ], 
 \qquad t\in \real_+ 
. 
\end{equation} 
\end{prop} 
\begin{Proof}
 By construction of $(X (t))_{t\in \real_+}$ 
 in Relations~\eqref{a1.x}-\eqref{a4.x} 
 we have 
$$ 
\P(X_k (t) = 1 \mid X_k (0) )
 = p_k\left(
 1 + e^{-t}Y_k \sqrt{\frac{q_k}{p_k}}\right)
,
$$ 
\\
$$ 
 \P(X_k (t) = -1 \mid X_k (0) )
 = q_k\left( 1-e^{-t}Y_k \sqrt{\frac{p_k}{q_k}}\right)
, 
$$ 
\\
$$ 
\P(X_k (t) = 1 \mid X_k (0) )
 = p_k\left(
 1 + e^{-t}Y_k \sqrt{\frac{q_k}{p_k}}\right)
,
$$ 
\\
$$ 
 \P(X_k (t) = -1 \mid X_k (0) )
 = q_k\left( 1-e^{-t}Y_k \sqrt{\frac{p_k}{q_k}}\right)
, 
$$ 
 thus 
$$ 
 d \P(X_k (t)(\tilde{\omega} ) = \epsilon \mid  X (0)) (\omega ) 
 = 
 \left( 1 + e^{-t}Y_k (\omega ) Y_k(\tilde{\omega} ) 
 \right) 
 d\P(X_k (\tilde{\omega} ) = \epsilon ) , 
$$ 
 $\varepsilon = \pm 1$. 
 Since the components of $(X_k (t))_{k\in\inte}$ are 
 independent, this shows that 
 the law of $(X_0 (t),\ldots ,X_n (t))$ conditionally
 to $X (0)$ has the density $q_t^n(\tilde{\omega},\cdot )$ 
 with respect to $\P$:
\begin{eqnarray*} 
\lefteqn{ 
\! \! \! \! \! \! \! \! \! \! \! \! \! \! \! 
 d\P(X_0 (t)(\tilde{\omega} ) = \epsilon_0, 
 \ldots, X_n (t)(\tilde{\omega} ) =\epsilon_n \mid X (0) ) (\tilde{\omega} ) 
} 
\\ 
 & = & q^n_t(\tilde{\omega},\omega )
 d\P(X_0 (\tilde{\omega} ) = \epsilon_0, \ldots, 
 X_n (\tilde{\omega} ) = \epsilon_n ). 
\end{eqnarray*} 
 Consequently we have 
\begin{equation} 
\label{rprf} 
 \ee[F ( X (t) ) \mid X (0) = \tilde{\omega} ] 
 = \int_\Omega F(\omega ) q^N_t (\tilde{\omega} , \omega )
 \P (d\omega) 
, 
\end{equation} 
 hence from \eqref{repr}, Relation~\eqref{istheou} holds 
 for $F\in L^2(\Omega , {\cal F}_N)$, $N\geq 0$. 
\end{Proof} 
 The independent components $X_k (t)$, $k\in \inte$, 
 can be constructed from the data of 
 $X_k (0) =\epsilon$ and an independent exponential random variable $\tau_k$ 
 via the following procedure. 
 If $\tau_k <t$, let $X_k (t) =X_k (0) =\epsilon$, 
 otherwise if $\tau_k >t$, take $X_k (t)$ to be an independent copy of 
 $X_k$. 
 This procedure is illustrated in the following equalities: 
\begin{eqnarray} 
\label{r1r} 
 \P(X_k (t) = 1 \mid X_k (0) = 1) 
 & = & \ee[ {\bf 1}_{\{\tau_k >t\}} ] 
 + \ee[ {\bf 1}_{\{\tau_k <t\}} {\bf 1}_{\{X_k=1\}}] 
\\ 
\nonumber 
 & = & e^{-t} + p_k (1-e^{-t}) 
, 
\\
\nonumber 
\\ 
\label{r2r} 
 \P(X_k (t) = -1 \mid X_k (0) = 1) 
 & = & \ee[ {\bf 1}_{\{\tau_k <t\}} {\bf 1}_{\{X_k=-1\}}] 
\\ 
\nonumber 
 & = & q_k (1-e^{-t}) , 
\\
\nonumber 
\\ 
\label{r3r} 
 \P(X_k (t) = -1 \mid X_k (0) = -1) 
 & = & \ee[ {\bf 1}_{\{\tau_k >t\}} ] 
 + \ee[ {\bf 1}_{\{\tau_k <t\}} {\bf 1}_{\{X_k=-1\}}] 
\\
\nonumber 
& = & e^{-t} + q_k (1-e^{-t}) , 
\\ 
\nonumber 
\\ 
\label{r4r} 
 \P(X_k (t) = 1 \mid X_k (0) = -1) 
 & = & \ee[ {\bf 1}_{\{\tau_k <t\}} {\bf 1}_{\{X_k=1\}}] 
\\ 
\nonumber 
 & = & p_k (1-e^{-t}) 
. 
\end{eqnarray} 
 The operator $L^2 (\Omega \times \inte ) \to 
 L^2 (\Omega \times \inte )$ 
 which maps $(u_k)_{k\in\inte}$ to $(P_tu_k)_{k\in\inte}$ 
 is also denoted by $P_t$. 
 As a consequence of the representation of $P_t$ given 
 in Lemma~\ref{p11} we obtain the following bound. 
\begin{lemma} 
\label{bbdd}
 For $F\in \Dom (D)$ we have 
$$
\Vert P_t u\Vert_{L^\infty (\Omega , \ell^2(\inte ))}
 \leq \Vert u\Vert_{L^\infty (\Omega , \ell^2(\inte ))},
 \qquad 
 t\in \real_+, 
 \quad 
 u\in L^2 (\Omega \times \inte ) 
. 
$$
\end{lemma}
\begin{Proof}
 As a consequence of the representation formula \eqref{rprf} we have 
 $\P(d \tilde{\omega} )$-a.s.: 
\begin{eqnarray*} 
\Vert P_t u \Vert_{\ell^2(\inte )}^2 ( \tilde{\omega} ) 
 & = & 
\sum_{k=0}^\infty
 | P_t u_k ( \tilde{\omega} ) |^2
\\ 
& = & \sum_{k=0}^\infty
 \left( \int_{\Omega} 
 u_k ( \omega ) Q_t ( \tilde{\omega} , d\omega ) \right)^2
\\ 
& \leq  &
 \sum_{k=0}^\infty
 \int_{\Omega} 
 | u_k ( \omega )|^2 Q_t ( \tilde{\omega} , d\omega ) 
\\ 
& =  &
 \int_{\Omega} 
 \Vert u\Vert_{\ell^2(\inte )}^2 ( \omega ) 
 Q_t ( \tilde{\omega} , d\omega ) 
\\ 
& \leq & \Vert u\Vert_{L^\infty (\Omega , \ell^2(\inte ))}^2 
. 
\end{eqnarray*} 
\end{Proof} 
\section{Covariance Identities}\index{covariance identities!discrete time}
\label{s8} 
 In this section we state the covariance identities 
 which will be used for the proof of deviation inequalities 
 in the next section. 
 The covariance 
 $\Cov (F,G)$ of $F,G\in L^2 (\Omega )$
 is defined as
$$\Cov (F,G) = \ee[(F-\ee[F])(G-\ee[G])] 
 = \ee[FG] - \ee[F]\ee[G] 
. 
$$ 
\begin{prop}
\label{covclark0} 
 We have for
 $F,G\in L^2 (\Omega )$ such that $\ee[\Vert DF \Vert_{\ell^2(\inte )}^2]< \infty$: 
\begin{equation}
\label{covclark}
\Cov (F,G) =
 \ee\left[ \sum_{k=0}^\infty
 \ee\left[ D_{k} G \mid {\cal F}_{k-1} \right] D_k F 
 \right].
\end{equation}
\end{prop}
\begin{Proof} 
 This identity is a consequence of the Clark formula \eqref{clk}: 
\begin{eqnarray*}
\Cov (F,G) & = & \ee[(F-\ee[F])(G-\ee[G])]
\\
& = &
\ee\left[
\left( 
\sum_{k=0}^\infty \ee[D_kF \mid {\cal F}_{k-1}] Y_k
\right) 
\left( 
\sum_{l=0}^\infty \ee[D_lG \mid {\cal F}_{l-1}] Y_l
\right) 
\right]
\\
& = &
\ee\left[
\sum_{k=0}^\infty \ee[D_kF \mid {\cal F}_{k-1}]
 \ee[D_kG\mid {\cal F}_{k-1}]
\right]
\\
& = &
\sum_{k=0}^\infty
 \ee\left[ 
 \ee[\ee[D_kG \mid {\cal F}_{k-1}]
 D_kF\mid {\cal F}_{k-1}]
\right]
\\
& = & 
 \ee\left[ \sum_{k=0}^\infty 
 \ee[D_kG \mid {\cal F}_{k-1}] D_kF 
 \right]. 
\end{eqnarray*}
\end{Proof} 
 A covariance identity can also be obtained
 using the semi-group $(P_t)_{t\in \real_+}$.
\begin{prop} 
\label{crsg} 
 For any $F,G\in L^2 (\Omega )$ such that 
$$\ee[\Vert DF\Vert_{\ell^2(\inte )}^2] < \infty 
 \qquad 
 \mbox{and} 
 \qquad 
 \ee[\Vert DG\Vert_{\ell^2(\inte )}^2] < \infty 
, 
$$ 
 we have 
\begin{equation} 
\label{coov} 
\Cov (F,G) = \ee\left[ 
 \sum_{k=0}^\infty
 \int_0^\infty
 e^{-t} ( D_k F ) P_t D_k G dt
 \right]. 
\end{equation} 
\end{prop}
\begin{Proof} 
 Consider $F = J_n (f_n)$ and $G= J_m(g_m )$. 
 We have
\begin{eqnarray*} 
\lefteqn{ 
\Cov (J_n(f_n),J_m(g_m)) 
 = 
\ee\left[
J_n(f_n)
J_m(g_m)
\right]
} 
\\
& = &
 {\bf 1}_{\{n=m\}}
n! \langle f_n
 , g_n {\bf 1}_{\Delta_n}
 \rangle_{\ell^2(\inte^n)}
\\ 
 & =& 
 {\bf 1}_{\{n=m\}}
 n! n 
 \int_0^\infty e^{-nt} dt
 \langle f_n
 , g_n {\bf 1}_{\Delta_n}
 \rangle_{\ell^2(\inte^n)} 
\\
& = & 
 {\bf 1}_{\{n-1=m-1\}}
 n! n  
 \int_0^\infty
e^{-t}
\sum_{k=0}^{\infty }
 \langle f_n (*,k)
 , e^{-(n-1)t}
 g_n (*,k)
 {\bf 1}_{\Delta_n} (*,k)
 \rangle_{\ell^2(\inte^{n-1})}
dt
\\
& = &
 n m 
 \ee\left[
\int_0^\infty
e^{- t} 
\sum_{k=0}^{\infty }
  J_{n-1}( f_n(*,k) {\bf 1}_{\Delta_n} (*,k) )
 e^{- ( m - 1 ) t} J_{m-1}( g_m(*,k) {\bf 1}_{\Delta_m} (*,k) )
dt
\right]
\\
& = &
 n m 
 \ee\left[
\int_0^\infty
 e^{-t} 
 \sum_{k=0}^{\infty }
 J_{n-1}( f_n(*,k) {\bf 1}_{\Delta_n}(*,k) )
 P_t J_{m-1}( g_m(*,k) {\bf 1}_{\Delta_m}(*,k) )
dt
\right]
\\ 
 & = & 
\ee\left[
\int_0^\infty
e^{-t}
\sum_{k=0}^{\infty }
D_k J_n(f_n) P_t D_k J_m(g_m)
dt
\right]
. 
\end{eqnarray*}
\end{Proof} 
\noindent 
 From \eqref{r1r}-\eqref{r4r} 
 the covariance identity \eqref{coov} shows that 
\begin{eqnarray} 
\nonumber 
\lefteqn{ 
\Cov (F,G) = \ee\left[ 
 \sum_{k=0}^\infty
 \int_0^\infty 
 e^{-t} D_k F P_t D_k G dt
 \right] 
} 
\\ 
\nonumber 
& = & 
 \ee\left[ 
 \int_0^1 
 \sum_{k=0}^\infty 
 D_k F P_{-\log \alpha} D_k G d\alpha 
 \right] 
\\ 
\nonumber 
 & = & \int_0^1 
 \int_{\Omega \times \Omega} 
 \sum_{k=0}^\infty 
 D_k F (\omega ) 
 D_k G ((\omega_i {\bf 1}_{\{\tau_i < -\log \alpha\}} 
 + \omega'_i {\bf 1}_{\{\tau_i < -\log \alpha\}} )_{i\in\inte} ) 
 d\alpha \P(d\omega ) \P(d\omega') 
\\ 
\nonumber 
& = & 
 \int_0^1 
 \int_{\Omega \times \Omega} 
 \sum_{k=0}^\infty 
 D_k F (\omega ) 
 D_k G ((\omega_i {\bf 1}_{\{\xi_i < \alpha\}} 
 + \omega'_i {\bf 1}_{\{\xi_i > \alpha\}})_{i\in\inte} ) 
 \P(d\omega ) \P(d\omega')  d\alpha , 
\\ 
& & \label{gghhh} 
\end{eqnarray} 
 where $(\xi_i)_{i\in\inte}$ is a family of i.i.d. random variables, 
 uniformly distributed on $[0,1]$. 
 Note that the marginals of 
 $(X_k, X_k {\bf 1}_{\{\xi_k < \alpha\}} 
 + X_k' {\bf 1}_{\{\xi_i > \alpha\}})$ 
 are identical when $X_k'$ is an independent copy of 
 $X_k$. 
 Let 
$$\phi_\alpha (s,t) = \ee[e^{is X_k}e^{it (X_k+{\bf 1}_{\{ 
 \xi_k< \alpha \}}) + it (X_k'+{\bf 1}_{\{ 
 \xi_k > \alpha \}})}]. 
$$ 
 Then we have the relation 
$$\phi_\alpha (s,t) = 
 \alpha \phi (s+t) + (1-\alpha ) \phi (s) \phi (t), 
 \quad \alpha\in [0,1]. 
$$ 
 Note that 
$$ 
 \Cov (e^{isX_k},e^{itX_k}) 
 = \phi_1(s,t) - \phi_0 (s,t) 
 = \int_0^1 \frac{d \phi_\alpha }{d\alpha} (s,t) d\alpha 
 = \phi(s+t) - \phi(s)\phi(t). 
$$ 
 Next we prove an iterated version of the covariance identity 
 in discrete time, which is an analog of a result proved 
 in \cite{houdre} for the Wiener and Poisson processes. 
\begin{theorem} 
\label{th1a}  
 Let $n\in \inte$ and $F,G\in L^2 (\Omega )$.
 We have 
\begin{eqnarray} 
\label{*a} 
\lefteqn{ 
\! \! \! \! \! \! \! \! 
 \Cov (F,G) = 
 \sum_{d=1}^{d=n} 
 (-1)^{d+1} \ee\left[ 
 \sum_{\{1\leq k_1< \cdots < k_d\}} 
 (D_{k_d}\cdots D_{k_1} F)
 (D_{k_d}\cdots D_{k_1} G) 
 \right] 
} 
\\ 
\nonumber 
 & & 
 + 
 (-1)^n 
 \ee\left[ 
 \sum_{\{1\leq k_1< \cdots < k_{n+1} \}} 
 ( D_{k_{n+1}}\cdots D_{k_1} F ) 
 \ee\left[ D_{k_{n+1}}\cdots D_{k_1} G\mid {\cal F}_{k_{n+1}-1} \right] 
 \right].  
\end{eqnarray}
\end{theorem} 
\begin{Proof} 
 Take $F=G$. 
 For $n=0$, (\ref{*a}) is a consequence of the Clark formula. 
 Let $n\geq 1$. 
 Applying Lemma~\ref{lemmaa} to 
 $D_{k_n}\cdots D_{k_1} F$ with $a=k_n$ and 
 $b=k_{n+1}$, 
 and summing on 
 $(k_1, \ldots ,k_n)\in \Delta_n$,  
 we obtain 
\begin{eqnarray*} 
\lefteqn{ 
\! \! \! \! \! \! \! \! \! \! \! \! \! \! \! \! \! \! \! \! 
\! \! \! \! \! \! \! \! \! 
 \ee\left[  
 \sum_{\{1\leq k_1< \cdots < k_n \}} 
 \left( 
 \ee[D_{k_n}\cdots D_{k_1} F \mid {\cal F}_{k_n-1} ]\right)^2
 \right] 
 = 
 \ee\left[ 
 \sum_{\{1\leq k_1< \cdots < k_n \}} \mid D_{k_n}\cdots D_{k_1} F |^2 
 \right] 
} 
\\ 
& & 
 - 
 \ee\left[ 
 \sum_{\{1\leq k_1< \cdots < k_{n+1} \}} 
 \left( 
 \ee\left[ D_{k_{n+1}}\cdots D_{k_1} F\mid {\cal F}_{k_{n+1}-1} 
 \right]\right)^2  
 \right], 
\end{eqnarray*} 
 which concludes the proof by induction and bilinearity. 
\end{Proof} 
 As a consequence of Theorem~\ref{th1a}, 
 letting $F=G$ we get the variance 
 inequality 
$$\sum_{k=1}^{2n} \frac{(-1)^{k+1}}{k!} 
 \ee \left[ 
 \Vert D^k F\Vert_{\ell^2(\Delta_k)}^2 
 \right] 
 \leq \Var (F) \leq 
 \sum_{k=1}^{2n-1} 
 \frac{(-1)^{k+1}}{k!} 
 \ee \left[ 
 \Vert D^k F\Vert_{\ell^2(\Delta_k)}^2 
 \right] 
, 
$$ 
 since 
\begin{eqnarray*} 
\lefteqn{ 
\ee\left[ 
 \sum_{\{1\leq k_1< \cdots < k_{n+1} \}} 
 ( D_{k_{n+1}}\cdots D_{k_1} F ) 
 \ee\left[ D_{k_{n+1}}\cdots D_{k_1} G\mid {\cal F}_{k_{n+1}-1} \right] 
 \right] 
} 
\\ 
& = & 
\ee\left[ 
 \sum_{\{1\leq k_1< \cdots < k_{n+1} \}} 
  \ee\left[ 
 ( D_{k_{n+1}}\cdots D_{k_1} F ) 
 \ee\left[ D_{k_{n+1}}\cdots D_{k_1} G\mid {\cal F}_{k_{n+1}-1} \right] 
 \mid {\cal F}_{k_{n+1}-1} \right] 
 \right] 
\\ 
& = & 
\ee\left[ 
 \sum_{\{1\leq k_1< \cdots < k_{n+1} \}} 
 ( 
 \ee\left[ D_{k_{n+1}}\cdots D_{k_1} G\mid {\cal F}_{k_{n+1}-1} \right] 
 )^2 
\right] 
\\ 
& \geq & 0 
, 
\end{eqnarray*} 
 see Relation~(2.15) in \cite{houdre} in continuous time. 
 In a similar way, another iterated covariance identity can be 
 obtained from Proposition~\ref{crsg}. 
\begin{corollary} 
 Let $n\in \inte$ and $F,G \in L^2(\Omega, {\cal F}_N)$. 
 We have 
\begin{eqnarray} 
\nonumber 
\Cov (F,G) & = & 
 \sum_{d=1}^{d=n} 
 (-1)^{d+1} \ee\left[ 
 \sum_{\{1\leq k_1< \cdots < k_d\leq N\}} 
 (D_{k_d}\cdots D_{k_1} F)
 (D_{k_d}\cdots D_{k_1} G) 
 \right] 
\\ 
\nonumber 
 & & 
 \! \! \! \! \! \! \! \! \! \! \! \! \! \! \! \! \! \! \! \! \! 
 \! \! \! \! \! \! \! \! \! \! \! \! \! \! \! \! \! \! \! \! \! \! 
 + 
 (-1)^n 
 \int_{\Omega \times \Omega} 
 \sum_{\{1\leq k_1< \cdots < k_{n+1} \leq N\}} 
 D_{k_{n+1}}\cdots D_{k_1} F(\omega ) 
 D_{k_{n+1}}\cdots D_{k_1} G(\omega') 
\\ 
\label{****a} 
 & & 
 q^N_t(\omega,\omega')\P(d\omega )\P(d\omega' ).  
\end{eqnarray} 
\end{corollary} 
 The covariance and variance 
 have the tensorization property: 
$$\Var (FG) = \ee[F\Var G] + \ee[G\Var F]$$ 
 if $F,G$ are independent, hence most of the identities in this 
 section can be obtained by tensorization of a one dimensional 
 elementary covariance identity. 
\bigskip 

 An elementary consequence of the covariance identities is 
 the following lemma. 
\begin{lemma} 
\label{m2} 
 Let $F, G\in L^2 (\Omega )$ such that 
$$\ee [D_k F | {\cal F}_{k-1} ] 
 \cdot 
 \ee [D_k G | {\cal F}_{k-1} ]\ge 0, 
 \qquad
 k\in\inte. 
$$ 
 Then $F$ and $G$ are non-negatively correlated: 
$$\Cov (F, G)\ge 0. 
$$ 
\end{lemma} 
\noindent 
 According to the next definition, a non-decreasing functional $F$ 
 satisfies $D_kF\ge 0$ for all $k\in\inte$. 
\begin{definition}  
  A random variable $F : \Omega \to \real$ is said to be non-decreasing if 
 for all $\omega_1,\omega_2\in \Omega$ we have 
$$\omega_1 (k) \le \omega_2 (k), \qquad 
 k\in\inte, 
 \qquad \Rightarrow 
 \ 
 F(\omega_1 ) \leq F(\omega_2) 
. 
$$ 
\end{definition} 
\noindent 
 The following result is then immediate from 
 Proposition~\ref{fimdtl} and Lemma~\ref{m2}, 
 and shows that the FKG inequality holds on $\Omega$. 
 It can also be obtained from from Proposition~\ref{crsg}.  
\begin{prop} 
 If $F, G\in L^2 (\Omega )$ are non-decreasing 
 then $F$ and $G$ are non-negatively correlated: 
$$ 
\Cov (F, G)\ge 0. 
$$ 
\end{prop} 
\noindent 
 Note however that the assumptions of 
 Lemma~\ref{m2} are actually weaker as they 
 do not require $F$ and $G$ to be non-decreasing. 
\section{Deviation Inequalities}\index{deviation inequalities!discrete time}
\label{devsec} 
 In this section, which is based on \cite{hp}, 
 we recover a deviation inequality of \cite{bobkov} in the case 
 of Bernoulli measures, using covariance representations instead 
 of the logarithmic  Sobolev inequalities to be presented 
 in Section~\ref{lsid}. 
 The method relies on a bound on the Laplace transform 
 $L(t) = \ee [ e^{t F} ]$ obtained via a differential inequality 
 and Chebychev's inequality. 
\begin{prop}
\label{c1.02}
 Let $F : \Omega \to \real$ be such that
 $\vert F_k^+ - F_k^- \vert \leq K$, $k\in\inte$, for
 some $K\geq 0$, and
 $\Vert DF \Vert_{L^\infty (\Omega , \ell^2(\inte ))} <\infty$.
 Then
\begin{eqnarray*}
\P(F-\ee[F]\geq x) & \leq & \exp \left(
 - \frac{\Vert DF\Vert_{L^\infty (\Omega ,\ell^2 (\inte ))}^2}{K^2}
 g\left( \frac{x K}{\Vert DF\Vert_{L^\infty (\Omega ,\ell^2 (\inte ))}^2}\right)
 \right)
\\
& \leq & \exp \left(
 - \frac{x}{2 K}
 \log \left(
 1 + \frac{x K}{\Vert DF\Vert_{L^\infty (\Omega ,\ell^2 (\inte ))}^2}
 \right)\right),
\end{eqnarray*}
 with $g(u) = (1+u)\log (1+u)-u$, $u\geq 0$.
\end{prop}
\begin{Proof} 
 Although $D_k$ does not satisfy a derivation rule for products, from 
 Proposition~\ref{chnrle} we have 
\begin{eqnarray*}
D_k e^F & = &
 {\bf 1}_{\{X_k=1\}}
 \sqrt{p_kq_k}
 (e^F- e^{F^-_k})
 + {\bf 1}_{\{X_k=-1\}}
 \sqrt{p_kq_k}
 (e^{F^+_k}- e^F)
\\
& = &
 {\bf 1}_{\{X_k=1\}}
 \sqrt{p_kq_k}
 e^F (1- e^{-\frac{1}{\sqrt{p_kq_k}}D_kF})
 + {\bf 1}_{\{X_k=-1\}}
  \sqrt{p_kq_k} e^F (e^{\frac{1}{\sqrt{p_kq_k}} D_kF}- 1)
\\
& = &
 - X_k \sqrt{p_kq_k}
 e^F (e^{-\frac{X_k}{\sqrt{p_kq_k}}D_kF}- 1) 
, 
\end{eqnarray*}
 hence 
\begin{equation}
\label{prd.0}
D_k e^F =
 X_k \sqrt{p_kq_k}e^F
 (1-e^{-\frac{X_k}{\sqrt{p_kq_k}} D_kF}), 
\end{equation}
 and  since the function $x \mapsto (e^x-1)/x$ is positive
 and increasing on $\real$ we have:
$$\frac{e^{-sF}D_ke^{sF}}{D_kF}
 = -\frac{X_k\sqrt{p_kq_k} }{D_kF}
 \left(
 e^{- s \frac{X_k}{\sqrt{p_kq_k}} D_kF} -1 \right)
 \leq
 \frac{e^{sK}-1}{K},$$ 
 or in other terms: 
$$ 
 \frac{e^{-sF}D_ke^{sF}}{D_kF} 
 = {\bf 1}_{\{X_k=1 \}} 
 \frac{e^{s ( F_k^--F_k^+) }-1}{F_k^--F_k^+} 
 + {\bf 1}_{\{X_k= -1 \}} 
 \frac{e^{s ( F_k^+-F_k^-) }-1}{F_k^+-F_k^-} 
 \leq
 \frac{e^{sK}-1}{K}.$$ 
 We first assume that $F$ is a bounded random variable with $\ee[F]=0$. 
 From Lemma~\ref{bbdd} applied to $DF$, we have
\begin{eqnarray*} 
\ee[Fe^{sF}] & =& \Cov (F,e^{sF}) 
\\ 
 & = & \ee\left[\int_0^\infty
 e^{-v}
 \sum_{k=0}^\infty D_k e^{sF}
 P_v D_k F dv \right]
\\
& \leq &
 \left\| \frac{e^{-sF}
 De^{sF}}{DF} \right\|_{\infty}
 \ee\left[e^{sF}
 \int_0^\infty e^{-v}
 \Vert DFP_v DF \Vert_{\ell^1(\inte )} dv
 \right]
\\ 
& \leq &
 \frac{e^{sK}-1}{K}
 \ee\left[e^{sF}
 \Vert DF \Vert_{\ell^2 (\inte )} 
 \int_0^\infty e^{-v}
 \Vert P_v DF \Vert_{\ell^2(\inte )} dv
 \right]
\\
& \leq &
 \frac{e^{sK}-1}{K}
 \ee\left[e^{sF} \right]
 \Vert DF \Vert_{L^\infty (\Omega , \ell^2(\inte ))}^2
 \int_0^\infty e^{-v} dv
\\
& \leq &
 \frac{e^{sK}-1}{K}
 \ee\left[e^{sF} \right]
 \Vert DF \Vert_{L^\infty(\Omega ,\ell^2(\inte ))}^2.
\end{eqnarray*} 
 In the general case, letting $L(s) = \ee[e^{s(F-\ee[F])}]$, we have 
\begin{eqnarray*} 
\log (\ee[e^{t(F-\ee[F])}])
 & = & 
 \int_0^t \frac{L'(s)}{ L(s)}
 ds
\\ 
& = & 
 \int_0^t
 \frac{\ee[(F-\ee[F]) e^{s(F-\ee[F])}]}{\ee[e^{s(F-\ee[F])}]}
 ds
\\ 
&  \leq & 
 \frac{1}{K} \Vert DF \Vert_{L^\infty (\Omega ,\ell^2 (\inte ))}^2 
 \int_0^t (e^{sK}-1) ds
\\ 
& = & 
 \frac{1}{K^2} (e^{tK}-tK-1)
 \Vert DF \Vert_{L^\infty (\Omega ,\ell^2 (\inte ))}^2 
,
\end{eqnarray*} 
$t\geq 0$.
 We have for all $x\geq 0$ and $t\geq 0$:
\begin{eqnarray*} 
\P(F-\ee[F]\geq x)
& \leq & e^{-tx} \ee[e^{t(F-\ee[F])}]
\\ 
& \leq & 
 \exp \left( 
 \frac{1}{K^2} (e^{tK}-tK-1)
 \Vert DF \Vert_{L^\infty (\Omega ,\ell^2 (\inte ))}^2 
 -tx \right) 
, 
\end{eqnarray*} 
 The minimum in $t\geq 0$ in the above expression is attained with 
$$t = \frac{1}{K} 
 \log \left( 1 + \frac{xK}{ \Vert DF \Vert_{L^\infty (\Omega ,\ell^2 (\inte ))}^2 }\right), 
$$ 
 hence 
\begin{eqnarray*} 
\lefteqn{ 
\P(F-\ee[F]\geq x)
} 
\\ 
& \leq & \exp \left( 
 -\frac{1}{K} \left(
 \left( 
 x+ \frac{1}{K} \Vert DF \Vert_{L^\infty (\Omega ,\ell^2 (\inte ))}^2 \right)
 \log \left( 1+ x K 
 \Vert DF \Vert_{L^\infty (\Omega ,\ell^2 (\inte ))}^{-2} \right)
 - x\right)
\right) 
\\
 & \leq & 
 \exp \left( 
 -\frac{x}{2K}
 \log \left( 1+ x K
 \Vert DF \Vert_{L^\infty (\Omega ,\ell^2 (\inte ))}^{-2}
 \right) \right) ,
\end{eqnarray*}
 where we used the inequality
 $(1+u)\log (1+u)-u \geq \frac{u}{2}\log (1+u)$.
 If $K=0$, the above proof is still valid by replacing all terms
 by their limits as $K\to 0$.
 If $F$ is not bounded the conclusion
 holds for $F_n = \max (-n , \min ( F, n))$, $n\geq 1$,
 and $(F_n)_{n\in\inte}$,
 $(DF_n)_{n\in\inte}$, converge respectively almost surely 
 and in $L^2(\Omega \times \inte )$ to $F$ and $DF$, with
 $\Vert DF_n \Vert_{L^\infty (\Omega ,L^2 (\inte ))}^2
 \leq \Vert DF \Vert_{L^\infty (\Omega ,L^2 (\inte ))}^2$.
\end{Proof}
 In case $p_k=p$ for all $k\in\inte$,
 the conditions
$$\frac{1}{\sqrt{pq}} \vert D_k F\vert \leq \beta, \quad 
 k\in\inte, 
 \quad 
 \mbox{and} 
 \quad 
 \Vert DF\Vert_{L^\infty (\Omega ,\ell^2 (\inte ))}^2 \leq \alpha^2 
,
$$ 
 give
$$
\P(F-\ee[F]\geq x) \leq \exp \left(
 - \frac{\alpha^2 pq}{\beta^2}
 g\left( \frac{x \beta}{\alpha^2 pq}\right)
 \right)
 \leq \exp \left(
 - \frac{x}{2 \beta}
 \log \left(
 1 + \frac{x \beta}{\alpha^2 pq}
 \right)\right),
$$
 which is Relation (13) in \cite{bobkov}.
 In particular if $F$ is ${\cal F}_N$-measurable, then
$$
\P(F-\ee[F]\geq x) \leq \exp \left(
 -Ng\left( \frac{x}{\beta N}\right)
\right)
 \leq \exp \left(
 -\frac{x}{\beta }
 \left(
 \log\left(
 1+\frac{x}{\beta N}\right)
 -1\right)
\right).
$$
 Finally we show a Gaussian concentration inequality
 for functionals of $(S_n)_{n\in\inte}$, using the covariance
 identity \eqref{covclark}. 
 We refer to \cite{bobkovsyk}, \cite{bht2}, \cite{ht2},
 \cite{ledouxesaim}, for other versions of this inequality.
\begin{prop}
\label{gsian}
 Let $F: \Omega \rightarrow \real$ be such that
$$
 \left\|
 \sum_{k=0}^\infty
 \frac{1}{2(p_k\wedge q_k)}
 \vert  D_k F \vert
 \Vert  D_k F \Vert_\infty
 \right\|_\infty
 \leq K^2.
$$
 Then
\begin{equation}
\label{bdd}
\P(F-\ee[F]\geq x) \leq \exp \left(
 - \frac{x^2}{2
 K^2 }
 \right), \qquad x\geq 0.
\end{equation}
\end{prop} 
\begin{Proof} 
 Again, we assume that $F$ is a bounded random variable with $\ee[F] = 0$.
 Using the inequality 
\begin{equation}
\label{1..}
\vert e^{tx}-e^{ty} \vert
 \leq \frac{t}{2} 
 \vert x-y \vert (e^{tx}+e^{ty}),
 \qquad x, y\in \real,
\end{equation}
 we have 
\begin{eqnarray}
  \nonumber
\vert D_k e^{tF} \vert & = &
 \sqrt{p_kq_k}
 \vert e^{tF_k^+}- e^{tF_k^-} \vert
\\ 
  \nonumber
& \leq & \frac{1}{2}
 \sqrt{p_kq_k}
 t \vert F_k^+ - F_k^- \vert
 (e^{tF_k^+} + e^{tF_k^-} )
\\
  \nonumber
& = & \frac{1}{2}
 t \vert D_k F \vert
 (e^{tF_k^+} + e^{tF_k^-} )
\\ 
\label{djkddsad}
 & \leq & \frac{t}{2(p_k\wedge q_k)}
 \vert D_k F \vert
\ee\left[
 e^{tF}
 \mid X_i, \ i\not=k
 \right]
\\
  \nonumber
 & = &
 \frac{1}{2(p_k\wedge q_k)}
 t \ee\left[
 e^{tF}
 \vert D_k F \vert
 \mid X_i, \ i\not=k
 \right], 
\end{eqnarray}
 where in \eqref{djkddsad}
 the inequality is due to the absence of chain
 rule of derivation for the operator $D_k$.
 Now, Proposition~\ref{covclark0} yields 
\begin{eqnarray*} 
\ee[Fe^{tF}] & = & \Cov (F,e^{s F}) 
\\ 
 &=& 
 \sum_{k=0}^\infty
 \ee[\ee[D_k F \mid {\cal F}_{k-1} ]
 D_k e^{tF} ]
\\ 
& \leq & \sum_{k=0}^\infty
 \Vert  D_k F \Vert_\infty
 \ee\left[\vert D_k e^{tF} \vert \right]
\\
 & \leq & \frac{t}{2} \sum_{k=0}^\infty
 \frac{1}{p_k\wedge q_k}
 \Vert  D_k F \Vert_\infty
 \ee\left[ 
 \ee\left[
 e^{tF}
 \vert D_k F \vert
 \mid X_i, \ i\not=k
 \right] \right]
\\
 & = & \frac{t}{2} 
  \ee\left[
  e^{tF}
 \sum_{k=0}^\infty 
 \frac{1}{p_k\wedge q_k}
 \Vert  D_k F \Vert_\infty
 \vert D_k F \vert
 \right]
\\
 & \leq & \frac{t}{2} 
 \ee[ e^{tF} ]
 \left\|
 \sum_{k=0}^\infty 
 \frac{1}{p_k\wedge q_k}
 \vert  D_k F \vert
 \Vert  D_k F \Vert_\infty
 \right\|_\infty.
\end{eqnarray*}
 This shows that
\begin{eqnarray*}
 \log (\ee[e^{t(F-\ee[F])}]) 
 & = & \int_0^t \frac{\ee[(F-\ee[F])e^{s(F-\ee[F])}]}{\ee[e^{s(F-\ee[F])}]} ds 
\\ 
& \leq & K^2 \int_0^t s ds 
\\ 
 & = & \frac{t^2}{2} K^2,
\end{eqnarray*}
 hence 
\begin{eqnarray*} 
 e^{x} \P(F-\ee[F]\geq x)
 & \leq & \ee[e^{t(F-\ee[F])}]
\\ 
& \leq & e^{t^2 K^2 / 2 }, \qquad t\geq 0,
\end{eqnarray*} 
 and
$$ 
 \P(F-\ee[F]\geq x) 
 \leq e^{\frac{t^2}{2} K^2-tx}, \quad t\geq 0.
$$ 
 The best inequality is obtained for $t=x/K^2$.
 If $F$ is not bounded the conclusion
 holds for $F_n = \max (-n , \min ( F, n))$, $n\geq 0$,
 and $(F_n)_{n\in\inte}$, 
 $(DF_n)_{n\in\inte}$, converge respectively to $F$ and $DF$
 in $L^2(\Omega )$, resp. $L^2(\Omega \times \inte )$, with
 $\Vert DF_n \Vert_{L^\infty (\Omega ,\ell^2 (\inte ))}^2
 \leq \Vert DF \Vert_{L^\infty (\Omega ,\ell^2 (\inte ))}^2$.
\end{Proof}
 The bound \eqref{bdd} implies
 $\ee[e^{\alpha \vert F\vert}]<\infty$ for all $\alpha >0$,
 and $\ee[e^{\alpha F^2}]<\infty$ for all $\alpha <
 1/(2K^2)$.
 In case $p_k=p$, $k\in\inte$, we obtain
$$\P(F-\ee[F]\geq x) \leq \exp \left(
 - \frac{px^2}{\Vert DF \Vert_{
 \ell^2(\inte,L^\infty (\Omega ))}^2}
 \right).
$$
\section{Logarithmic Sobolev Inequalities}\index{logarithmic Sobolev inequalities!discrete time} 
\label{lsid} 
 The logarithmic Sobolev inequalities on Gaussian space 
 provide an infinite dimensional analog of Sobolev inequalities, 
 cf. e.g. \cite{ledouxmarkov}. 
 On Riemannian path space \cite{capitaine} and on Poisson space \cite{ane}, 
 \cite{wuls2}, martingale methods have been successfully applied to 
 the proof of logarithmic Sobolev inequalities. 
 Here, discrete time martingale methods are used 
 along with the Clark predictable representation formula 
 \eqref{clk} as in \cite{gaopri}, 
 to provide a proof of logarithmic Sobolev inequalities 
 for Bernoulli measures. 
 Here we are only concerned with modified logarithmic 
 Sobolev inequalities, and we refer to 
 \cite{saloff}, Theorem~2.2.8 and references 
 therein, for the standard version of the 
 logarithmic Sobolev inequality on the hypercube 
 under Bernoulli measures. 
\bigskip 

\begin{figure}[!ht]
\begin{center} 
\resizebox*{10cm}{7cm}{\rotatebox{0}{\includegraphics[scale=0.9]{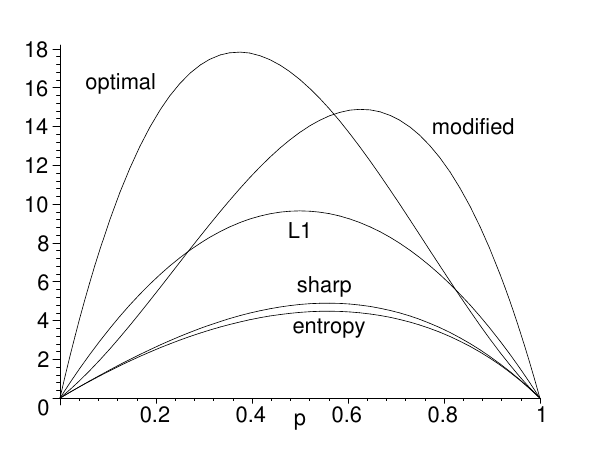}}} 
\end{center} 
\centering
\caption{\small 
 Graph of the entropy as a function of $p$. 
}
\label{g5.0} 
\end{figure}

 The entropy of a random variable 
 $F>0$ is defined by
$$\Ent [F] = \ee [F \log F ] -
 \ee [F] \log \ee [F],$$
 for sufficiently integrable $F$. 
\begin{lemma} 
\label{tensorlemma} 
 The entropy has the tensorization 
 property, i.e. if $F,G$ are sufficiently integrable 
 independent random variables we have 
\begin{equation} 
\label{tensop} 
\Ent [ FG ] 
 = 
\ee[F\Ent [ G ] ] + \ee[G\Ent [F] ] 
. 
\end{equation} 
\end{lemma}
\begin{Proof} 
 We have 
\begin{eqnarray*} 
\nonumber 
\Ent [ FG ] & = & \ee[FG\log (FG)] - \ee[FG]\log \ee[FG] 
\\ 
\nonumber 
& = & \ee[FG ( \log F + \log G) ] - \ee[F]\ee[G] ( \log \ee[F] + \log \ee[G] ) 
\\ 
\nonumber 
& = & \ee[G] \ee[F \log F ] + \ee[F] \ee[ G \log G) ] - \ee[F]\ee[G] ( \log \ee[F] + \log \ee[G] ) 
\\ 
 & = & 
\ee[F\Ent [ G ] ] + \ee[G\Ent [F] ] 
. 
\end{eqnarray*} 
\end{Proof} 
 In the next proposition we recover the 
 modified logarithmic Sobolev inequality of \cite{bobkov} 
 using the Clark representation formula in discrete time. 
\begin{theorem}
\label{lsaprop}
 Let $F\in \Dom (D)$
 with $F>\eta$ a.s. for some $\eta>0$. 
 We have 
\begin{equation}
\label{lsa}
\Ent [F]
 \leq \ee
 \left[\frac{1}{F} \Vert DF\Vert_{\ell^2(\inte )}^2
 \right].
\end{equation}
\end{theorem}
\begin{Proof}
 Assume that $F$ is ${\cal F}_N$-measurable and 
 let $M_n = \ee [F\mid {\cal F}_n]$, $0\leq n\leq N$.
 Using Corollary~\ref{cormes} and the Clark formula \eqref{clk} 
 we have 
$$M_n = M_{-1} + \sum_{k=0}^n 
 u_k Y_k, 
 \qquad 
 0 \leq n \leq N 
, 
$$
 with $u_k = \ee [D_k F \mid {\cal F}_{k-1}]$,
 $0\leq k \leq n \leq N$,
 and $M_{-1}=\ee[F]$.
 Letting $f(x) = x\log x$ and using the bound
\begin{eqnarray*} 
 f(x+y)-f(x) & = & y\log x + (x+y)\log \left(
 1+\frac{y}{x} \right)
\\ 
& \leq & y(1+\log x) + \frac{y^2}{x},
\end{eqnarray*} 
 we have:
\begin{eqnarray*}
 \Ent [F]
 & = &
  \ee[f(M_N)]-\ee[f(M_{-1})]
\\
& = &
 \ee \left[
\sum_{k=0}^N 
f(M_k )
-f(M_{k-1})
\right]
\\
&  = &
\ee \left[
\sum_{k=0}^N 
f\left(
M_{k-1} + Y_k u_k
\right)
-f(M_{k-1})
\right]
\\
& \leq &
 \ee \left[
\sum_{k=0}^N 
 Y_k u_k (1+\log M_{k-1} ) + \frac{Y^2_k u^2_k}{M_{k-1}}
\right]
\\
 & = &
 \ee \left[
\sum_{k=0}^N 
\frac{1}{\ee[F\mid {\cal F}_{k-1}]}
 (\ee[D_kF\mid {\cal F}_{k-1}])^2
\right]
\\
 & \leq &
 \ee \left[
\sum_{k=0}^N 
 \ee\left[\frac{1}{F}
 | D_kF |^2 \mid {\cal F}_{k-1}\right]
\right]
\\ 
 & = & 
 \ee \left[
\frac{1}{F}
 \sum_{k=0}^N 
 | D_kF |^2
\right]. 
\end{eqnarray*}
 where we used the Jensen inequality 
 and the convexity of $(u,v)\mapsto v^2/u$ on $(0,\infty)\times \real$, 
 or the Schwarz inequality applied to $1/\sqrt{F}$ 
 and $(D_k F/\sqrt{F})_{k\in \inte}$, 
 as in the Wiener and Poisson cases \cite{capitaine} and \cite{ane}.
 This inequality is extended by density to $F\in \Dom (D)$. 
\end{Proof} 
 Theorem~\ref{lsaprop} can also be recovered by 
 the tensorization Lemma~\ref{tensorlemma} and the following 
 one-variable argument: 
 letting $p+q=1$, $p,q>0$, $f:\{ -1,1\} \to (0,\infty )$,
 $\ee[f]=pf(1)+qf(-1)$, and $\d f = f(1)-f(-1)$ we have:
\begin{eqnarray*}
\lefteqn{ 
\Ent [ f ] = pf(1)\log f(1) + qf(-1)\log f(-1)
 - \ee[f]\log \ee[f]
}
\\
& = & pf(1)\log (\ee[f] + q\d f)
 + qf(-1)\log (\ee[f] -p \d f)
 - (pf(1)+qf(-1))\log \ee[f]
\\
& = & pf(1)\log \left( 1 + q \frac{\d f}{\ee[f]}\right)
 + qf(-1)\log \left( 1 - p \frac{\d f}{\ee[f]}\right)
\\
& \leq & pqf(1) \frac{\d f}{\ee[f]}
 - pqf(-1) \frac{\d f}{\ee[f]}
 = pq \frac{ | \d f |^2}{\ee[f]}
\\
& \leq & pq \ee\left[ \frac{1}{f} | \d f |^2
 \right]
. 
\end{eqnarray*}
 Similarly we have 
\begin{eqnarray*}
\Ent [ f ] & = & pf(1)\log f(1) + qf(-1)\log f(-1)
 - \ee[f]\log \ee[f]
\\
& = & p (\ee[f] + q\d f) \log (\ee[f] + q\d f)
\\
& & + q(\ee[f] -p \d f) \log (\ee[f] -p \d f)
 - (pf(1)+qf(-1))\log \ee[f]
\\
& = & p \ee[f] \log \left(1 + q\frac{\d f}{\ee[f]}\right)
 + p q\d f \log f(1)
\\
& & + q\ee[f] \log \left(1 -p \frac{\d f}{\ee[f]}\right)
 - qp \d f \log f(-1)
\\
& \leq &
 p q\d f \log f(1)
 - pq \d f \log f(-1)
\\
& = & pq \ee\left[ \d f \d \log f \right]
, 
\end{eqnarray*} 
 which, by tensorization, recovers the following $L^1$ inequality of 
 \cite{gao}, \cite{daipra}, and proved in \cite{wuls2} in the Poisson 
 case. 
 In the next proposition we state and prove 
 this inequality in the multidimensional case, 
 using the Clark representation formula, similarly 
 to Theorem~\ref{lsaprop}. 
\begin{theorem}
\label{lsal1} 
 Let $F>0$ be ${\cal F}_N$-measurable.
 We have 
\begin{equation}
\label{lsa3}
\Ent [F]
 \leq \ee
 \left[\sum_{k=0}^N 
 D_kF D_k \log F \right]
.
\end{equation}
\end{theorem}
\begin{Proof} 
 Let $f(x) = x\log x$ and 
$$\Psi (x,y) = (x+y)\log (x+y) - x\log x - (1+\log x)y,
 \quad x, \ x+y >0. 
$$ 
 From the relation 
\begin{eqnarray*} 
\lefteqn{ 
Y_ku_k = Y_k \ee [ D_k F \mid {\cal F}_{k-1} ] 
} 
\\ 
& = & q_k {\bf 1}_{\{X_k=1\}} 
 \ee[(F_k^+-F_k^-)\mid {\cal F}_{k-1}] 
 + p_k {\bf 1}_{\{X_k=-1\}} 
 \ee[(F_k^- - F_k^+)\mid {\cal F}_{k-1}] 
\\
& = & 
{\bf 1}_{\{X_k=1\}} 
\ee[(F_k^+-F_k^-){\bf 1}_{\{X_k=-1\}} 
 \mid {\cal F}_{k-1} ] 
+ {\bf 1}_{\{X_k= -1\}} 
 \ee[(F_k^--F_k^+){\bf 1}_{\{X_k=1\}} 
 \mid {\cal F}_{k-1} ],  
\end{eqnarray*} 
 we have, using the convexity of $\Psi$: 
\begin{eqnarray*}
\lefteqn{ \Ent [F] = \ee \left[
\sum_{k=0}^N 
f\left(
M_{k-1} + Y_k u_k
\right)
-f(M_{k-1})
\right]
} 
\\
& = &
 \ee \left[
\sum_{k=0}^N 
\Psi (
M_{k-1}, Y_k u_k
)
+ Y_k u_k (1+\log M_{k-1} )
\right]
\\
& = &
 \ee \left[
\sum_{k=0}^N 
\Psi (
M_{k-1}, Y_k u_k
)
\right]
\\
& = &
 \ee \left[
 \sum_{k=0}^N 
 p_k
 \Psi
 \left(
 \ee[F\mid {\cal F}_{k-1}] 
 ,
 \ee[(F_k^+-F_k^-){\bf 1}_{\{X_k=-1\}} 
 \mid {\cal F}_{k-1} ] 
 \right)
\right. 
\\
& & \left. 
 +
 q_k
 \Psi \left( \ee[F\mid {\cal F}_{k-1}] 
 ,
 \ee[(F_k^--F_k^+){\bf 1}_{\{X_k=1\}} \mid 
 {\cal F}_{k-1} ] 
 \right)
 \right]
\\ 
& \leq & 
 \ee \left[
 \sum_{k=0}^N 
 \ee\left[ 
 p_k
 \Psi
 \left(
 F 
 ,
 (F_k^+-F_k^-){\bf 1}_{\{X_k=-1\}} 
 \right)
 +
 q_k
 \Psi \left( F 
 ,
 (F_k^--F_k^+){\bf 1}_{\{X_k=1\}} 
 \right)
 \mid {\cal F}_{k-1} 
 \right]\right]
\\
& = &
 \ee \left[
 \sum_{k=0}^N 
 p_k{\bf 1}_{\{X_k=-1\}} 
 \Psi
 \left(
 F_k^- 
 ,
 F_k^+-F_k^-
 \right)
 +
 q_k{\bf 1}_{\{X_k=1\}} 
 \Psi \left( F_k^+ 
 ,
 F_k^--F_k^+
 \right)
 \right]
\\
& = &
 \ee \left[
 \sum_{k=0}^N 
 p_k q_k
 \Psi(F_k^-
 ,
 F_k^+ - F_k^-
 )
 +
 p_kq_k
 \Psi(F_k^+
 ,
 F_k^- - F_k^+
 )
 \right]
\\
& = &
 \ee \left[
 \sum_{k=0}^N 
 p_k q_k
 (\log F_k^+ - \log F_k^- )
 (F_k^+ - F_k^- )
\right]
\\
& = &
 \ee \left[
 \sum_{k=0}^N 
 D_k F D_k \log F
\right]
. 
\end{eqnarray*}
\end{Proof} 
 The proof of Theorem~\ref{lsal1} can also be obtained 
 by first using the bound
$$f(x+y)-f(x) = y\log x + (x+y)\log \left(
 1+\frac{y}{x} \right)
 \leq y(1+\log x) + y\log (x+y),
$$ 
 and then the convexity of $(u,v)\to v(\log (u+v)-\log u)$:
\begin{eqnarray*} 
\lefteqn{
 \Ent [F]
 =
 \ee \left[
\sum_{k=0}^N 
f\left(
M_{k-1} + Y_k u_k
\right)
-f(M_{k-1})
\right]
}
\\
& \leq &
 \ee \left[
\sum_{k=0}^N 
 Y_k u_k (1+\log M_{k-1} ) +
 Y_k u_k
 \log (
 M_{k-1} + Y_k u_k )
\right]
\\
 & = &
 \ee \left[
\sum_{k=0}^N 
 Y_k u_k
 (\log (
 M_{k-1} + Y_k u_k )
 -
 \log M_{k-1} )
\right]
\\ 
 & = &
 \ee \left[
\sum_{k=0}^N 
 \sqrt{p_kq_k} 
 \ee[D_kF\mid {\cal F}_{k-1}] 
 (\log 
 \ee[F + (F_k^+-F_k^-){\bf 1}_{\{X_k=-1\}} 
 \mid {\cal F}_{k-1} ] 
 -
 \log \ee[F\mid {\cal F}_{k-1}] )
\right. 
\\ 
& & \left. 
 - \sqrt{p_kq_k} 
 \ee[D_kF\mid {\cal F}_{k-1}] 
 (\log 
 \ee[F +(F_k^--F_k^+){\bf 1}_{\{X_k=-1\}} \mid {\cal F}_{k-1} ] 
 -
 \log \ee[F\mid {\cal F}_{k-1} ])\right]
\\
 & \leq & 
 \ee \left[ 
\sum_{k=0}^N 
 \ee \left[ 
 \sqrt{p_kq_k} 
 D_kF 
 (\log (
 F + (F_k^+-F_k^-){\bf 1}_{\{X_k=-1\}} 
 ) -
 \log F )
\right. \right. 
\\ 
& & \left. \left. 
 - \sqrt{p_kq_k} 
 D_kF 
 (\log (
 F +(F_k^--F_k^+){\bf 1}_{\{X_k=1\}} 
 )
 -
 \log F)
 \mid {\cal F}_{k-1} 
\right]\right]
\\ 
 & = &
 \ee \left[
\sum_{k=0}^N 
 \sqrt{p_kq_k} 
 D_kF {\bf 1}_{\{X_k=-1\}} 
 (\log F_k^+ -
 \log F_k^- )
 \right. 
\\ 
 & & 
 \left. 
 - \sqrt{p_kq_k} 
 D_kF {\bf 1}_{\{X_k=1\}} 
 (\log F_k^- - \log F_k^+ )
\right]
\\
 & \leq &
 \ee \left[
\sum_{k=0}^N 
 \sqrt{p_kq_k} 
 q_k D_kF 
 (\log F_k^+ 
 - \log F_k^- )
 - \sqrt{p_kq_k} 
 p_k D_kF 
 (\log F_k^- - \log F_k^+ )
\right]
\\
& = &
 \ee \left[
\sum_{k=0}^N 
 D_k F D_k \log F
 \right]
. 
\end{eqnarray*}
 The application of Theorem~\ref{lsal1} to $e^F$ gives 
 the following inequality for $F>0$, ${\cal F}_N$-measurable: 
\begin{eqnarray} 
\nonumber 
\Ent[e^F] & \leq &  \ee \left[
 \sum_{k=0}^N 
 D_k F D_k e^F
\right]
\\ 
\nonumber 
& = &
 \ee \left[
 \sum_{k=0}^N 
 p_k q_k
 \Psi(e^{F_k^-}
 ,
 e^{F_k^+} - e^{F_k^-}
 )
 +
 p_k q_k
 \Psi( e^{F_k^+}
 ,
 e^{F_k^-} - e^{F_k^+}
 )
 \right]
\\
\nonumber 
& = &
 \ee \left[
 \sum_{k=0}^N 
 p_k q_k 
 e^{F_k^-}
 ((F_k^+-F_k^-)e^{F_k^+-F_k^-}
 - e^{F_k^+-F_k^-}
 + 1 )
\right.
\\
\nonumber 
 & & \left.  +
 p_k q_k 
 e^{F_k^+}
 ((F_k^- -F_k^+ )e^{F_k^- -F_k^+}
 - e^{F_k^--F_k^+}
 + 1 ) \right]
\\
\nonumber 
& = &
 \ee \left[
 \sum_{k=0}^N 
 p_k {\bf 1}_{\{X_k=-1\}}
 e^{F_k^-}
 ((F_k^+-F_k^-)e^{F_k^+-F_k^-}
 - e^{F_k^+-F_k^-}
 + 1 )
\right.
\\
\nonumber 
 & & \left.  +
 q_k {\bf 1}_{\{X_k=1\}}
 e^{F_k^+}
 ((F_k^- -F_k^+ )e^{F_k^- -F_k^+}
 - e^{F_k^--F_k^+}
 + 1 ) \right]
\\
\label{lsa3.0} 
& = &
 \ee \left[
 e^{F}
 \sum_{k=0}^N 
 \sqrt{p_kq_k}
 \vert Y_k \vert
 ( \nabla_k F e^{ \nabla_k F }
 - e^{\nabla_k F} +1)
 \right]
. 
\end{eqnarray}
 This implies
\begin{equation}
\label{lll}
 \Ent [e^F]
 \leq
 \ee \left[
 e^{F}
 \sum_{k=0}^N 
 ( \nabla_k F e^{ \nabla_k F }
 - e^{\nabla_k F} +1)
 \right]. 
\end{equation} 
\noindent 
 As already noted in \cite{daipra}, 
 \eqref{lsa3} and the Poisson limit theorem 
 yield the $L^1$ inequality of \cite{wuls2}.
 Let $M_n = (n + X_1 +\cdots +X_n)/2$,
 $F = \varphi ( M_n )$,
 and $p_k = \lambda /n$, $k\in \inte$, $\lambda >0$.
 Then
\begin{eqnarray*}
\lefteqn{ 
\sum_{k=0}^n 
 D_kF D_k \log F
} 
\\ 
 & = & \frac{\lambda}{n}
 \left(1-\frac{\lambda}{n} \right)
 (n-M_n) (\varphi (M_n + 1) - \varphi (M_n))
 \log (\varphi (M_n + 1) - \varphi (M_n))
\\
& & +
 \frac{\lambda}{n}\left(1-\frac{\lambda}{n}\right)
 M_n (\varphi (M_n ) - \varphi (M_n - 1))
 \log (\varphi (M_n ) - \varphi (M_n - 1)).
\end{eqnarray*}
 In the limit we obtain
$$\Ent [\varphi (U) ] \leq \lambda \ee
 [ (\varphi (U+1)-\varphi (U))(\log \varphi (U+1)-\log \varphi (U))],
$$
 where $U$ is a Poisson random variable with parameter $\lambda$.
 In one variable we have, still letting $\d f = f(1)-f(-1)$, 
\begin{eqnarray*} 
\Ent [ e^f ] & \leq & pq \ee\left[ \d e^f \d \log e^f \right] 
\\ 
& = & p q(e^{f(1)}-e^{f(-1)}) (f(1)-f(-1))
\\
& = & p qe^{f(-1)}(
 (f(1)-f(-1))e^{f(1)-f(-1)}
 - e^{f(1)-f(-1)}
 + 1 )
\\
& & +
 p qe^{f(1)}(
 (f(-1)-f(1))e^{f(-1)-f(1)}
 - e^{f(-1)-f(1)}
 + 1 )
\\
& \leq & qe^{f(-1)}(
 (f(1)-f(-1))e^{f(1)-f(-1)}
 - e^{f(1)-f(-1)}
 + 1 )
\\
& & +
 p e^{f(1)}(
 (f(-1)-f(1))e^{f(-1)-f(1)}
 - e^{f(-1)-f(1)}
 + 1 )
\\
& = & \ee[e^f(
 \nabla f e^{\nabla f}
 - e^{\nabla f}
 + 1 )
]
, 
\end{eqnarray*} 
 where $\nabla_k$ is the gradient operator defined in 
 \eqref{mod2}. 
 This last inequality is not comparable to the optimal constant inequality
\begin{equation}
\label{sh}
\Ent [e^F]
 \leq
 \ee \left[
 e^{F}
 \sum_{k=0}^N 
 p_kq_k ( \vert \nabla_k F \vert e^{\vert \nabla_k F \vert }
 - e^{\vert \nabla_k F \vert} +1)
 \right]
,
\end{equation}
 of \cite{bobkov} since when $F_k^+-F_k^-\geq 0$ the right-hand side of
 \eqref{sh} grows as $F_k^+e^{2F_k^+}$,
 instead of $F_k^+e^{F_k^+}$ in \eqref{lll}. 
 In fact we can prove the following inequality which improves 
 \eqref{lsa}, \eqref{lsa3} and \eqref{sh}. 
\begin{theorem} 
\label{lsal1.1} 
 Let $F$ be ${\cal F}_N$-measurable. We have
\begin{equation} 
\label{lsa3.01}
 \Ent [e^F]
 \leq
 \ee \left[ e^{F}
 \sum_{k=0}^N 
 p_kq_k 
 ( \nabla_k F e^{ \nabla_k F }
 - e^{\nabla_k F} +1)
 \right] 
. 
\end{equation}
\end{theorem} 
 Clearly, \eqref{lsa3.01} is better than 
 \eqref{sh}, \eqref{lsa3.0} and \eqref{lsa3}. 
 It also improves \eqref{lsa} from the bound 
$$xe^x-e^x+1\leq ( e^x-1 )^2, \quad x\in \real 
, 
$$ 
 which implies 
$$e^F ( \nabla F e^{\nabla F} -e^{\nabla F} + 1 ) \leq e^F ( e^{\nabla F}-1 )^2 
 = 
 e^{-F} | \nabla e^{F} |^2
. 
$$ 
 By the tensorization property \eqref{tensop}, the proof of 
 \eqref{lsa3.01} reduces to the following 
 one dimensional lemma. 
\begin{lemma} For any $0\leq p\leq 1$, $t\in{\mathbb R}$, $a \in{\mathbb
R}$, $q=1-p$, 
\begin{eqnarray*} 
\lefteqn{ 
pt e^t+qa e^a-\left(p e^t+q e^a\right)\log\left( p
e^t+ q e^a\right)
}
\\
&\leq & p q \left( q e^a\left((t-a)e^{t-a}-e^{t-a}+1\right)
+pe^t\left((a-t)e^{a-t}-e^{a-t}+1\right)\right).
\end{eqnarray*} 
\end{lemma}
\begin{Proof} Set 
\begin{eqnarray*} 
g(t) &= &pq \left(qe^a\left((t-a)e^{t-a}-e^{t-a}+1\right)
+pe^t\left((a-t)e^{a-t}-e^{a-t}+1\right)\right)\\
& & - pt e^t-qa e^a+\left(p e^t+qe^a\right)\log\left( p
e^t+qe^a\right).
\end{eqnarray*} 
 Then
$$g'(t)= pq\left(q e^a(t - a)e^{t - a} + pe^t\left(-e^{a-t}
+ 1\right)\right) 
 - pte^t+pe^t\log(pe^t + q e^a)
$$ 
and $g''(t)=p e^t h(t)$, 
where
$$ 
h(t) = - a - 2pt - p + 2pa + p^2t - p^2a 
 + \log( pe^t + q e^a ) + \frac{pe^t}{p e^t + q e^a }.
$$ 
 Now, 
\begin{eqnarray*} 
h'(t) & =&- 2p + p^2+ \frac{ 2 pe^t}{pe^t + qe^a} -
\frac{p^2e^{2t}}{(p e^t + qe^a )^2}\\
&=&\frac{pq^2(e^t-e^a)(pe^t+(q+1)e^a )}{(pe^t+qe^a)^2},
\end{eqnarray*} 
which implies that $h'(a)=0$, $h'(t)<0 $ for any $t<a$ and
$h'(t)>0 $ for any $t>a$. 
 Hence, for any $t\not= a$,
$h(t)>h(a)=0$, and so $g''(t)\geq 0$ for any $t\in {\mathbb R}$
and $g''(t)=0$ if and only if $t=a$. Therefore, $g'$ is
strictly increasing. Finally, since $t=a$ is the unique root of
$g'=0$, we have that $g(t)\geq g(a)=0$ for all $t\in \real$. 
\end{Proof} 
\noindent 
 This inequality improves \eqref{lsa}, \eqref{lsa3}, and \eqref{sh}, 
 as illustrated in one dimension in Figure~\ref{g5.0}, where 
 the entropy is represented as a function of $p\in [0,1]$ 
 with $f(1)=1$ and $f(-1)=3.5$. 
 The inequality \eqref{lsa3.01} is a discrete analog of the sharp inequality 
 on Poisson space of \cite{wuls2}. 
 In the symmetric case $p_k=q_k=1/2$, $k\in \inte$, 
 we have 
\begin{eqnarray*} 
 \Ent [e^F]
 & \leq & 
 \ee\left[ 
 e^F 
 \sum_{k=0}^N 
 p_kq_k 
 ( \nabla_k F e^{ \nabla_k F }
 -  \nabla_k F +1)
\right] 
\\ 
 & = & 
 \frac{1}{8} 
 \ee \left[
 \sum_{k=0}^N 
 e^{F_k^-}
 ((F_k^+-F_k^-)e^{F_k^+-F_k^-}
 - e^{F_k^+-F_k^-}
 + 1 )
\right.
\\
 & & \left.  +
  e^{F_k^+}
 ((F_k^- -F_k^+ )e^{F_k^- -F_k^+}
 - e^{F_k^--F_k^+}
 + 1 ) \right]
\\ 
 & = & 
 \frac{1}{8} 
 \ee \left[
 \sum_{k=0}^N 
 (e^{F_k^+}-e^{F_k^-}) 
 (F_k^+ -F_k^- ) 
\right]
\\
& = & 
 \frac{1}{2} 
 \ee\left[ 
 \sum_{k=0}^N 
 D_k F D_k e^ F 
 \right], 
\end{eqnarray*} 
 which improves on \eqref{lsa3}. 
\\ 

 Letting $F=\varphi (M_n)$ we have 
\begin{eqnarray*} 
\lefteqn{ 
 \ee\left[ 
 e^F 
 \sum_{k=0}^N 
 p_kq_k 
 ( \nabla_k F e^{ \nabla_k F }
 -  \nabla_k F +1)
 \right] 
} 
\\
& = &
 \frac{\lambda}{n}
 \left(1-\frac{\lambda}{n} \right)
 \ee\left[ 
 M_n
 e^{\varphi (M_n )}
\right. 
\\ 
 & & 
 \times 
 \left. 
 ((\varphi (M_n )
 - \varphi (M_n - 1)
 )
 e^{\varphi (M_n )
 - \varphi (M_n - 1)
 }
 - e^{\varphi (M_n )
 - \varphi (M_n - 1)
 }
 + 1 ) 
 \right] 
\\
 & & 
 \! \! \! \! \! \! \! \! \! \! \! \! \! \! 
 +
 \frac{\lambda}{n}
 \left(1-\frac{\lambda}{n} \right)
 \ee\left[ 
 (n-M_n) e^{\varphi (M_n )}
 \right. 
\\ 
 & & 
 \times 
 \left. 
 ((\varphi (M_n +1) - \varphi (M_n )
 ) 
 e^{\varphi (M_n +1 )
 - \varphi (M_n )}
 - e^{\varphi (M_n +1 )
 - \varphi (M_n )
 } + 1 ) 
\right] 
, 
\end{eqnarray*} 
 and in the limit as $n$ goes to infinity we obtain
$$\Ent [ e^{\varphi (U)} ] \leq \lambda \ee
 [ e^{\varphi (U)} (
 (\varphi (U+1)-\varphi (U)) e^{\varphi (U+1)-\varphi (U)}
 - e^{\varphi (U+1)-\varphi (U)}
 +1
 )
],
$$
 where $U$ is a Poisson random variable with parameter $\lambda$.
 This corresponds to the sharp inequality of \cite{wuls2}.
\section{Change of Variable Formula} 
\label{s11} 
 In this section we state a discrete-time analog of It\^o's change 
 of variable formula which will be useful for the predictable 
 representation of random variables and for option hedging. 
\begin{prop} 
\label{it} 
 Let $(M_n)_{n\in \inte}$ be a square-integrable martingale 
 and $f : \real\times \inte \to \real$. 
 We have 
\begin{eqnarray} 
\label{note00} 
\lefteqn{ 
 f(M_n,n) 
} 
\\ 
\nonumber 
 & =  & 
f(M_{-1},-1)
 + \sum_{k=0}^n 
 D_k f(M_k,k) 
 Y_k 
 + 
\sum_{k=0}^n 
 \ee[f(M_k,k) -f(M_{k-1},k-1) \mid {\cal F}_{k-1}] 
. 
\end{eqnarray} 
\end{prop}
\begin{Proof}
 By Proposition~\ref{martrepr} there exists square-integrable 
 process $(u_k)_{k\in\inte}$ such that 
$$M_n = M_{-1} + \sum_{k=0}^n u_k Y_k, 
 \qquad 
 n\in\inte. 
$$ 
 We write 
\begin{eqnarray*}
 f(M_n,n) - f(M_{-1},-1) & = & \sum_{k=0}^n 
 f(M_k,k)-f(M_{k-1},k-1) 
\\
& = &
\sum_{k=0}^n 
 f(M_k,k)-f(M_{k-1},k) + f(M_{k-1},k)-f(M_{k-1},k-1)
\\
& = &
 \sum_{k=0}^n 
 \sqrt{\frac{p_k}{q_k}}
\left(
f\left(
M_{k-1} + u_k \sqrt{\frac{q_k}{p_k}}, k
\right)
-f(M_{k-1},k)
\right)
Y_k
\\
& &
+
\frac{p_k}{q_k}
{\bf 1}_{\{X_k=- 1 \}}
\left(
f\left(
M_{k-1} + u_k \sqrt{\frac{q_k}{p_k}},k
\right)
-f(M_{k-1},k)
\right)
\\
& & +
{\bf 1}_{\{X_k=- 1 \}}
\left(
f\left(
M_{k-1} - u_k \sqrt{\frac{p_k}{q_k}}
, k
\right)
-f(M_{k-1},k)
\right)
\\ 
 & & + 
 \sum_{k=0}^n 
 f(M_{k-1},k)-f(M_{k-1},k-1)
\\ 
 & = &
 \sum_{k=0}^n 
\sqrt{\frac{p_k}{q_k}}
\left(
f\left(
M_{k-1} + u_k \sqrt{\frac{q_k}{p_k}},k
\right)
-f(M_{k-1},k)
\right)
Y_k
\\
\nonumber
& & +
\sum_{k=0}^n 
\frac{1}{q_k}
{\bf 1}_{\{ X_k=-1 \}}
\ee[f(M_k,k)-f(M_{k-1},k) \mid {\cal F}_{k-1}]
\\
\nonumber
& & +
\sum_{k=0}^n 
f(M_{k-1},k)-f(M_{k-1},k-1)
. 
\end{eqnarray*}
 Similarly we have
\begin{eqnarray*}
f(M_n,n) & = &
f(M_{-1},-1)
- \sum_{k=0}^n 
\sqrt{\frac{q_k}{p_k}}
\left(
f\left(
M_{k-1} - u_k \sqrt{\frac{p_k}{q_k}}
,k 
\right)
-f(M_{k-1} 
,
k 
) 
\right)
Y_k
\\
& &
+ 
\sum_{k=0}^n 
\frac{1}{p_k}
{\bf 1}_{\{X_k= 1 \}}
\ee[f(M_k,k)-f(M_{k-1},k) \mid {\cal F}_{k-1}]
\\ 
 & & + 
 \sum_{k=0}^n 
 f(M_{k-1},k)-f(M_{k-1},k-1)
, 
\end{eqnarray*}
 Multiplying each increment in the above formulas 
 respectively by $q_k$ and $p_k$ and summing on $k$ we get 
\begin{eqnarray*}
\lefteqn{ 
f(M_n,n) = 
f(M_{-1},-1)
} 
\\ 
&  &+ \sum_{k=0}^n 
\sqrt{p_kq_k} 
\left(
f\left( 
_{k-1} + u_k \sqrt{\frac{q_k}{p_k}}
,k 
\right)
- 
f\left(
M_{k-1} - u_k \sqrt{\frac{p_k}{q_k}} 
,k 
\right)
\right)
Y_k
\\
& &
+ 
\sum_{k=0}^n 
\ee[f(M_k,k) \mid {\cal F}_{k-1}]
 -f(M_{k-1},k )
. 
\end{eqnarray*}
\end{Proof} 
 Note that in \eqref{note00} we have 
$$ D_k f(M_k,k) = 
\sqrt{p_kq_k} 
\left(
f\left( 
M_{k-1} + u_k \sqrt{\frac{q_k}{p_k}}
,k 
\right)
- 
f\left(
M_{k-1} - u_k \sqrt{\frac{p_k}{q_k}} 
,k 
\right)
\right)
, \quad k\in \inte. 
$$ 

\noindent 
 On the other hand, the term 
$$ 
\ee[f(M_k,k) -f(M_{k-1},k-1) \mid {\cal F}_{k-1}] 
$$ 
 is analog to the generator part in the continuous time 
 It\^o formula, and can be written as 
$$ 
 p_k 
f\left( 
M_{k-1} + u_k \sqrt{\frac{q_k}{p_k}}
,k 
\right)
 + 
 q_k 
f\left(
M_{k-1} - u_k \sqrt{\frac{p_k}{q_k}} 
,k 
\right)
 - 
f\left(
M_{k-1} ,k-1 
\right)
. 
$$ 
 When $p_n=q_n=1/2$, $n\in\inte$, we have 
\begin{eqnarray*} 
 f(M_n,n) & = & 
f(M_{-1},-1)
+ \sum_{k=0}^n 
\frac{f\left( 
M_{k-1} + u_k 
,k 
\right)
- 
f\left(
M_{k-1} - u_k 
,k 
\right)
} 
{2} 
Y_k 
\\ 
 & & 
 + 
\sum_{k=0}^n 
\frac{f\left( 
M_{k-1} + u_k 
,k 
\right)
+ 
f\left(
M_{k-1} - u_k 
,k 
\right)
 - 2 f\left(
M_{k-1} 
,k-1 
\right)
} 
{2} 
. 
\end{eqnarray*} 
 The above proposition also provides an explicit version of 
 the Doob decomposition for supermartingales. 
 Naturally if $(f(M_n,n))_{n\in\inte}$ is a martingale we have 
\begin{eqnarray*}
\lefteqn{ 
f(M_n,n) = 
f(M_{-1},-1)
} 
\\ 
 & & + \sum_{k=0}^n 
\sqrt{p_kq_k}
\left(
f\left(
M_{k-1} + u_k \sqrt{\frac{q_k}{p_k}}
, k
\right)
-
f\left(
M_{k-1} - u_k \sqrt{\frac{p_k}{q_k}}
, 
k 
\right)
\right)
Y_k
\\ 
 & = & f(M_{-1},-1) 
 + \sum_{k=0}^n 
 D_k f(M_k,k) Y_k 
. 
\end{eqnarray*}
 In this case the Clark formula, the martingale representation 
 formula Proposition~\ref{martrepr} and the change of variable formula all coincide. 
 In this case, we have in particular 
$$D_k f(M_k,k) 
 = \ee[D_k f(M_n,n) \mid {\cal F}_{k-1}] 
 = \ee[D_k f(M_k,k) \mid {\cal F}_{k-1}] 
, \qquad k\in \inte 
. 
$$ 
 If $F$ is an ${\cal F}_N$-measurable 
 random variable and $f$ is a function such that 
$$ \ee[F \mid {\cal F}_n] = f(M_n,n), \qquad -1 \leq n \leq N, 
$$ 
 we have $F = f(M_N,N)$, $\ee[F]= f(M_{-1},-1)$ and 
\begin{eqnarray*} 
 F & = & \ee[F] 
 + \sum_{k=0}^n 
 \ee[D_k f(M_N,N) \mid {\cal F}_{k-1}] 
 Y_k 
\\ 
& = & \ee[F] 
 + \sum_{k=0}^n 
 D_k f(M_k,k) 
 Y_k 
\\ 
& = & \ee[F] 
 + \sum_{k=0}^n 
 D_k \ee[f(M_N,N) \mid {\cal F}_k ] 
 Y_k 
. 
\end{eqnarray*} 
 Such a function $f$ exists if $(M_n)_{n\in\inte}$ is Markov and 
 $F=h(M_N)$. 
 In this case, consider the semi-group $(P_{k,n})_{0\leq k < n \leq N}$ 
 associated to $(M_n)_{n\in\inte}$ and defined by 
$$[P_{k,n} h](x) = \ee[h(M_n) \mid M_k = x] 
. 
$$ 
 Letting $f(x,n) = [P_{n,N}h](x)$ we can write 
$$ 
 F = \ee[F] 
 + \sum_{k=0}^n 
 \ee[D_k h(M_N) \mid {\cal F}_{k-1}] 
 Y_k 
 = \ee[F] 
 + \sum_{k=0}^n 
 D_k [P_{k,N}h (M_k) ] Y_k 
. 
$$ 
\section{Option Hedging in Discrete Time}\index{option hedging!discrete time} 
\label{hdg} 
 In this section we give a presentation of the Black-Scholes formula 
 in discrete time, or Cox-Ross-Rubinstein model, see e.g. 
 \cite{follmerschied}, 
 \cite{lamberton}, $\S$15-1 of \cite{williams}, or \cite{ruiz}, 
 as an application of the Clark formula. 
\\ 

 In order to be consistent with the notation of the previous 
 sections we choose to use the time scale $\inte$, hence the 
 index $0$ is that of the first random value of any stochastic 
 process, while the index $-1$ corresponds to its deterministic 
 initial value. 
\\    

 Let $(A_k)_{k\in \inte}$ be a riskless asset with initial value 
 $A_{-1}$, and defined by 
$$ 
A_n = A_{-1} \prod_{k=0}^n  (1+r_k), 
\qquad n\in \inte, 
$$ 
 where $(r_k)_{k\in \inte}$, is a sequence 
 of deterministic numbers such that $r_k>-1$, $k\in\inte$. 
 Consider a stock price with initial value $S_{-1}$, given 
 in discrete time as 
$$ 
 S_n = \left\{ 
 \begin{array}{ll} 
(1+b_n)S_{n-1}, & X_n=1, 
\\ 
\\
(1+a_n)S_{n-1}, & X_n=-1, \quad n\in \inte 
, 
\end{array} 
\right. 
$$ 
 where $(a_k)_{k\in \inte}$ and $(b_k)_{k\in \inte}$ are sequences 
 of deterministic numbers such that 
$$ 
 -1 < a_k < r_k < b_k, 
 \qquad 
 k\in \inte. 
$$ 
 We have 
$$S_n = S_{-1} \prod_{k=0}^n  
 \sqrt{(1+b_k)(1+a_k)} 
 \left( 
 \frac{1+b_k}{1+a_k}\right)^{X_k/2}, 
\quad n\in \inte 
. 
$$ 
 Consider now the discounted stock price given as 
\begin{eqnarray*} 
 \tilde{S}_n & = & 
 S_n \prod_{k=0}^n  
 (1+r_k)^{-1}  
\\ 
& = & 
 S_{-1} 
 \prod_{k=0}^n  
 \left( 
 \frac{1}{1+r_k} 
  \sqrt{(1+b_k)(1+a_k)} 
 \left( 
 \frac{1+b_k}{1+a_k}\right)^{X_k/2} 
 \right) 
, 
\quad n\in \inte 
. 
\end{eqnarray*} 
 If $-1 < a_k < r_k < b_k$, $k\in \inte$, then 
 $(\tilde{S}_n)_{n\in\inte}$ is a martingale with respect to 
 $({\cal F}_n)_{n\geq -1}$ under the probability $\P^*$ given by 
$$p_k = (r_k-a_k)/(b_k-a_k), 
 \quad 
 q_k = (b_k-r_k)/(b_k-a_k), \quad 
 k\in\inte. 
$$ 
 In other terms, under $\P^*$ we have 
$$ 
 \ee^* [ S_{n+1} \mid {\cal F}_n ] 
 = (1+r_{n+1}) S_n, 
 \qquad 
 n\geq -1 
, 
$$ 
 where $\ee^*$ denotes the expectation under 
 $\P^*$. 
 Recall that under this probability measure there is 
 absence of arbitrage and the market is complete. 
 From the change of variable formula Proposition~\ref{it} 
 or from the Clark formula \eqref{clk} 
 we have the martingale representation 
$$\tilde{S}_n = S_{-1} + \sum_{k=0}^n  
 Y_k D_k \tilde{S}_k 
 = S_{-1} + \sum_{k=0}^n  
 \tilde{S}_{k-1} 
 \sqrt{p_kq_k}\frac{b_k-a_k}{1+r_k} 
 Y_k 
. 
$$ 
\begin{definition} 
 A portfolio strategy is a pair of predictable processes 
 $(\eta_k)_{k\in\inte}$ and $(\zeta_k)_{k\in\inte}$ 
 where $\eta_k$, resp. $\zeta_k$ 
 represents the numbers of units invested over the 
 time period $(k,k+1]$ in the asset $S_k$, resp. $A_k$,  
 with $k\geq 0$. 
\end{definition} 
 The value at time $k \geq -1$ 
 of the portfolio $(\eta_k,\zeta_k)_{0\leq k \leq N}$ 
 is defined as 
\begin{equation}
\label{e43.0}
 V_k = \zeta_{k+1} A_k + \eta_{k+1} S_k, \qquad 
 k \geq -1
, 
\end{equation} 
 and its discounted value is defined as 
\begin{equation} 
\label{plm} 
\tilde{V}_n = V_n \displaystyle\prod_{k=0}^n (1+r_k)^{-1}, 
 \qquad 
 n\geq -1 
. 
\end{equation} 
\begin{definition} 
 A portfolio $(\eta_k,\zeta_k)_{k \in \inte}$ is said to be 
 self-financing if 
$$A_n(\zeta_{n+1}-\zeta_{n})+ S_n(\eta_{n+1}-\eta_{n} ) =0, 
 \qquad n\geq 0. 
$$ 
\end{definition} 
 Note that the self-financing condition implies 
$$ 
 V_n= \zeta_n A_n + \eta_n S_n, \qquad n\geq 0 
. 
$$ 
 Our goal is to hedge an arbitrary claim on $\Omega$, 
 i.e. given an ${\cal F}_N$-measurable 
 random variable $F$ we search for a portfolio 
 $(\eta_k,\zeta_k)_{ 0 \leq k \leq n}$ such that the equality 
\begin{equation} 
\label{prblm1} 
F = V_N = \zeta_N A_N + \eta_N S_N 
\end{equation} 
 holds at time $N\in \inte$. 
\begin{prop} 
 Assume that the portfolio 
 $(\eta_k,\zeta_k)_{0\leq k \leq N}$ is self-financing. 
 Then we have the decomposition 
\begin{equation} 
\label{e440.00}
 V_n = V_{-1}  \prod_{k=0}^n (1+r_k) 
 + \sum_{i=0}^n 
 \eta_i S_{i-1} 
 \sqrt{p_iq_i} (b_i-a_i) Y_i 
 \prod_{k=i+1}^n (1+r_k) 
. 
\end{equation}
\end{prop} 
\begin{Proof} 
 Under the self-financing assumption we have 
\begin{eqnarray*} 
\label{cll.0}
 V_i-V_{i-1} & = & \zeta_i (A_i-A_{i-1}) + \eta_i (S_i-S_{i-1}) 
\\ 
 & = & r_i \zeta_i A_{i-1} + (a_i{\bf 1}_{\{X_i=-1\}} + b_i{\bf 1}_{\{X_i=1\}}) \eta_i S_{i-1} 
\\ 
 & = & \eta_i S_{i-1} (a_i{\bf 1}_{\{X_i= -1\}} + b_i {\bf 1}_{\{X_i=1\}} -r_i) + r_i V_{i-1}
\\ 
 & = & \eta_i S_{i-1} \sqrt{p_iq_i} (b_i-a_i) Y_i + r_i V_{i-1} 
, \qquad i\in \inte, 
\end{eqnarray*} 
 hence for the discounted portfolio we get: 
\begin{eqnarray*} 
 \tilde{V}_i-\tilde{V}_{i-1} 
  & = & 
 \displaystyle\prod_{k=1}^i (1+r_k)^{-1} 
 {V}_i 
 - 
 \displaystyle\prod_{k=1}^{i-1} (1+r_k)^{-1} 
 {V}_{i-1} 
\\ 
  & = & 
 \displaystyle\prod_{k=1}^{i} (1+r_k)^{-1} 
 ( 
 {V}_i 
  - 
 {V}_{i-1} 
 - r_i {V}_{i-1} 
 ) 
\\ 
 & = &  \eta_i S_{i-1} 
 \sqrt{p_iq_i} (b_i-a_i) Y_i 
 \displaystyle\prod_{k=1}^i (1+r_k)^{-1} 
, \qquad i\in \inte, 
\end{eqnarray*} 
 which successively yields \eqref{e440.1.1} and \eqref{e440.00}. 
\end{Proof} 
 As a consequence of \eqref{e440.00} and \eqref{plm} 
 we immediately obtain 
\begin{equation} 
\label{e440.1.1} 
 \tilde{V}_n = \tilde{V}_{-1} 
 + \sum_{i=0}^n 
 \eta_i S_{i-1} 
 \sqrt{p_iq_i} (b_i-a_i) Y_i 
 \displaystyle\prod_{k=0}^i (1+r_k)^{-1} 
, 
 \qquad n\geq -1 
. 
\end{equation} 
 The next proposition provides a solution to the hedging 
 problem under the constraint \eqref{prblm1}. 
\begin{prop} 
 Given $F\in L^2(\Omega , {\cal F}_N )$, let 
\begin{equation} 
\label{eta} 
 \eta_n = \frac{1}{S_{n-1} \sqrt{p_nq_n} (b_n-a_n) } \ee^* 
 [D_n F \mid {\cal F}_{n-1} ]
 \prod_{k=n+1}^N (1+r_k)^{-1} 
, \qquad 0\leq n \leq N 
, 
\end{equation} 
 and 
\begin{equation} 
\label{zeta} 
\zeta_n = 
 A_n^{-1} 
 \left( 
 \displaystyle\prod_{k=n+1}^N (1+r_k)^{-1} 
 \ee^* [F \mid {\cal F}_n ] 
 -\eta_n S_n\right) , 
 \qquad 0\leq n \leq N 
. 
\end{equation} 
 Then the portfolio $( \eta_k , \zeta_k )_{0\leq k \leq n}$ 
 is self financing and satisfies 
$$ 
 \zeta_{n} A_n + \eta_{n} S_n 
 = 
 \prod_{k=n+1}^N (1+r_k)^{-1} 
 \ee^* [F \mid {\cal F}_n ] 
, \qquad 0 \leq n \leq N 
, 
$$ 
 in particular we have $V_N = F$, 
 hence $(\eta_k ,\zeta_k )_{0 \leq k \leq N}$ 
 is a hedging strategy leading to $F$. 
\end{prop} 
\begin{Proof} 
 Let $(\eta_k)_{-1\leq k \leq N}$ be defined by \eqref{eta} and 
 $\eta_{-1}=0$, and consider the process $(\zeta_n)_{0\leq n \leq N}$ defined by 
$$\zeta_{-1} = \frac{\ee^* [F ]}{S_{-1}} \prod_{k=0}^N (1+r_k)^{-1} 
 \quad 
 \mbox{and} 
 \quad
 \zeta_{k+1} = \zeta_k - \frac{(\eta_{k+1}-\eta_k)S_k}{A_k}, 
 \qquad 
 k=-1,\ldots , N-1 
. 
$$ 
 Then $(\eta_k , \zeta_k )_{-1\leq k \leq N}$ satisfies 
 the self-financing condition 
$$A_k(\zeta_{k+1}-\zeta_k)+ S_k(\eta_{k+1}-\eta_k ) =0, 
 \qquad -1 \leq k \leq N-1 
. 
$$ 
 Let now 
$$V_{-1} = 
 \ee^* [F ] 
 \prod_{k=0}^N (1+r_k)^{-1} 
, 
 \quad 
 \mbox{and} 
 \quad 
 V_n = \zeta_{n } A_n + \eta_n S_n, \qquad 0 \leq n \leq N, 
$$ 
 and 
$$\tilde{V}_n = V_n \displaystyle\prod_{k=0}^n (1+r_k)^{-1}, 
 \qquad -1\leq n \leq N 
. 
$$ 
 Since $(\eta_k,\zeta_k)_{-1 \leq k \leq N}$ is self-financing, 
 Relation~\eqref{e440.1.1} shows that 
\begin{equation} 
\label{tv} 
 \tilde{V}_n = 
 \tilde{V}_{-1} 
 + \sum_{i=0}^n 
 Y_i 
 \eta_i S_{i-1} 
 \sqrt{p_iq_i} (b_i-a_i) 
 \prod_{k=1}^i (1+r_k)^{-1} 
, 
 \qquad 
 -1\leq n \leq N   
. 
\end{equation} 
 On the other hand, from the Clark formula \eqref{clk} and the 
 definition of $(\eta_k)_{-1\leq k \leq N}$ we have 
\begin{eqnarray*} 
\lefteqn{ 
 \ee^* [ F \mid {\cal F}_n ] 
 \prod_{k=0}^N (1+r_k)^{-1} 
} 
\\ 
 & = & 
 \ee^* \left[ 
 \ee^* [F ] \prod_{k=0}^N (1+r_k)^{-1} 
 + 
 \sum_{i=0}^N 
 Y_i 
 \ee^* [D_i F \mid {\cal F}_{i-1} ] 
 \prod_{k=0}^N (1+r_k)^{-1} 
 \Big| 
 {\cal F}_n 
 \right] 
\\ 
 & = & 
 \ee^* [F ] \prod_{k=0}^N (1+r_k)^{-1} 
 + 
 \sum_{i=0}^n 
 Y_i 
 \ee^* [D_i F \mid {\cal F}_{i-1} ] 
 \prod_{k=0}^N (1+r_k)^{-1} 
\\ 
 & = & 
 \ee^* [F ] \prod_{k=0}^N (1+r_k)^{-1} 
 + \sum_{i=0}^n 
 Y_i 
 \eta_i S_{i-1} 
 \sqrt{p_iq_i} (b_i-a_i) 
 \prod_{k=1}^i (1+r_k)^{-1} 
\\ 
 & = & 
 \tilde{V}_n 
\end{eqnarray*} 
 from \eqref{tv}. 
 Hence 
$$ 
 \tilde{V}_n 
 = 
 \ee^* [ F \mid {\cal F}_n ] 
 \prod_{k=0}^N (1+r_k)^{-1} 
, 
 \qquad -1 \leq n \leq N, 
$$ 
 and 
$$ 
 V_{n} = 
 \ee^* [ F \mid {\cal F}_n ] 
 \prod_{k=n+1}^N (1+r_k)^{-1} 
, \qquad -1 \leq n \leq N. 
$$ 
 In particular we have $V_N = F$. 
 To conclude the proof we note that from the 
 relation $V_n = \zeta_n A_n + \eta_n S_n$, $0\leq n \leq N$, 
 the process $(\zeta_n)_{0\leq n \leq N}$ coincides 
 with 
 $(\zeta_n)_{0\leq n \leq N}$ defined 
 by \eqref{zeta}. 
\end{Proof} 
\noindent 
 Note that we also have 
$$ 
 \zeta_{n+1} A_n + \eta_{n+1} S_n 
 = 
 \ee^* [F \mid {\cal F}_n ] 
 \prod_{k=n+1}^N (1+r_k)^{-1} 
, \qquad -1 \leq n \leq N 
. 
$$ 
 The above proposition shows that there always exists a hedging strategy starting from 
$$\tilde{V}_{-1} = \ee^* [F] \prod_{k=0}^N (1+r_k)^{-1}. 
$$ 
 Conversely, if there exists a hedging strategy leading to 
$$\tilde{V}_N = F \prod_{k=0}^N (1+r_k)^{-1},$$ 
 then $(\tilde{V}_n)_{-1\leq n \leq N}$ is necessarily a martingale with 
 initial value 
$$\tilde{V}_{-1} = \ee^* [\tilde{V}_N ] = \ee^* [F] \prod_{k=0}^N (1+r_k)^{-1}. 
$$ 
 When $F=h(\tilde{S}_N)$, we have 
 $\ee^* [h(\tilde{S}_N) \mid {\cal F}_{k}] = f(\tilde{S}_k,k)$ with 
$$f(x,k ) = \ee^* \left[ 
 h\left( x 
 \prod_{i=k+1}^n 
 \frac{\sqrt{(1+b_k)(1+a_k)}}{1+r_k}  
 \left( 
 \frac{1+b_k}{1+a_k}\right)^{X_k/2} 
 \right) 
 \right] 
. 
$$ 
 The hedging strategy is given by 
\begin{eqnarray*} 
 \eta_k & = & 
 \frac{1}{S_{k-1} \sqrt{p_kq_k} (b_k-a_k) } 
 D_k f(\tilde{S}_k,k) 
 \prod_{i=k+1}^N (1+r_i)^{-1} 
\\ 
 & = & 
 \frac{\prod_{i=k+1}^N (1+r_i)^{-1} }{S_{k-1} (b_k-a_k) } 
 \left( 
 f\left(\tilde{S}_{k-1} 
  \frac{1+b_k}{1+r_k} , k\right) 
 - 
 f\left(\tilde{S}_{k-1} 
 \frac{1+a_k}{1+r_k} 
 , k\right) 
 \right) 
, \qquad k\geq -1. 
\end{eqnarray*} 
 Note that $\eta_k$ is non-negative (i.e. there is no short-selling) 
 when $f$ is an increasing function, e.g. in the case of European options 
 we have $f(x) = (x-K)^+$. 

 \footnotesize

 \def\cprime{$'$} \def\polhk#1{\setbox0=\hbox{#1}{\ooalign{\hidewidth
  \lower1.5ex\hbox{`}\hidewidth\crcr\unhbox0}}}
  \def\polhk#1{\setbox0=\hbox{#1}{\ooalign{\hidewidth
  \lower1.5ex\hbox{`}\hidewidth\crcr\unhbox0}}} \def\cprime{$'$}

\end{document}